%% file: main.tex
\documentclass{article}
\usepackage[a4paper, total={16cm, 23cm}]{geometry}

\input{ex_shared}

\begin{document}

\maketitle

\begin{abstract}
  This work concerns surrogate modeling for quantities of interest (\QoI)
  arising from parametric non-linear partial differential equations (PDEs).
  More specifically, we consider extending the sparse-grids surrogate modeling
  approach to incorporate derivatives of the \QoI with respect to the PDE
  parameters.
  We discuss why this operation is not straightforward and 
  propose a hybrid approach in which a sparse-grid scheme provides
  the collocation points in the parameter domain and  
  a suitable polynomial space, but the surrogate model is built with a least-squares approach. 
  We showcase our approach on several numerical tests, and we discuss in particular how its performance crucially
  depends on the relative cost and accuracy of evaluating the derivatives of the \QoI compared to evaluating the \QoI itself.
\end{abstract}

\section{Introduction}\label{sec:introduction}
A central task in computational science and engineering is to efficiently estimate
how the quantities of interest (\QoI) derived from the solution 
of a partial differential equation (PDE)
change upon varying the values of the PDE parameters.
This task is indeed at the core of \emph{many-query} applications such as uncertainty quantification, inverse problems, optimization, and real-time control. Since directly solving the PDE at hand for each value of the parameters is typically very computationally demanding, 
many \emph{surrogate modeling} techniques, that build 
approximations of a \QoI using only a limited number of computationally expensive
\QoI evaluations, have been proposed in literature to replace direct solvers
\cite{ghanem:UQbook}.

In this work, we consider the well-established sparse-grids surrogate modeling technique
\cite{babuska.nobile.eal:stochastic2,bungartz2004sparse,xiu.hesthaven:high},
that recasts the problem of building a surrogate model
for a \QoI as an interpolation problem over the parameter domain.
More specifically, our aim is to extend this methodology to incorporate
the values of the derivatives of the \QoI with respect to the PDE
parameters. This additional information can be effectively obtained 
for only a modest increase in computational cost,
once a value of the \QoI needed to build the standard (sparse-grid)
surrogate model has been computed.

Strategies to incorporate gradient information have been proposed in the last two decades for several surrogate modeling approaches,
from kernel methods such as Gaussian process regression and radial basis functions to polynomial chaos expansions, see \eg
\cite{bhaduri2020,chung2002, guo.eal:2018, jakeman2015, laurent2017, li2011, peng2016},
and more recently has been used also in the context of operator learning
\cite{gonzalez.sieiro.eal:deeponet,oleary-dino2024};
nonetheless, they are a relatively unexplored topic in the within the sparse-grids literature,
see \cite{deBaar2015,georg2020,jakeman.eal:adjoint-sg15}.
In particular, 
\cite{georg2020,jakeman.eal:adjoint-sg15} do not actually directly incorporate gradient information in the sparse-grid construction,
but rather employ derivative values as a sensitivity measure to drive the construction of a classical \emph{dimension-adaptive sparse grid}
\cite{gerstner.griebel:adaptive}.
Thus, \cite{deBaar2015} is the only work in literature that truly attempts to incorporate gradient information
in the sparse-grid construction, and it is also the starting point of our discussion.

It actually turns out that the approach proposed in \cite{deBaar2015} is unfortunately faulty; however,
the problem is not immediately evident in the reported numerical results,
since the devised construction still shows an apparent convergence.
We provide evidence of why the formulation in \cite{deBaar2015} cannot work,
and propose to solve the issue by moving to a hybrid approach, in which
a sparse-grid scheme provides the collocation points in the parameter domain and a polynomial space,
but the surrogate model is built with a least-squares approach.

The overall convenience of this modified technique fundamentally depends on
how the accuracy and cost of computing the partial derivatives of the \QoI
compare to those of computing the \QoI itself.
The significance of these factors is discussed and empirically examined using several numerical examples of increasing complexity:
a curated collection of closed-form functions and two versions
of a non-linear Darcy problem with uncertain permeability field.

This paper is structured as follows.
\Cref{sec:parametricPDEs} introduces the mathematical setting of parametric PDEs,
and discusses how the derivatives of a \QoI with respect to the parameters can be computed.
\Cref{sec:hermite} reviews key concepts of interpolation in both univariate and multivariate settings,
and discusses how to extend such concepts to incorporate gradient information.
\Cref{sec:SG} 
frames the approach of \cite{deBaar2015} as a nested application of the Smolyak
algorithm \cite{wolfers.tempone:smolyak},
shows that the problem actually lacks the necessary structure to apply this strategy,
and describes the new hybrid approach.
\Cref{sec:results}
showcases the proposed approach on several numerical tests, and compares its
effectiveness against standard sparse-grids.
Finally, \cref{sec:conclusions} concludes the paper.

\subsection{Parameter-dependent PDEs}\label{sec:parametricPDEs} 
Let $\sol\in\solSpace$ be the solution of a parametric PDE
\begin{equation}\label{eq:PDE}
\diffOp(\parVar, \sol) = 0,
\end{equation}
where $\diffOp$ is a non-linear differential operator, and the parameters $\parVar$ can take values in a Cartesian product of closed intervals
\[
  \parSpace=\bigotimes_{n=1}^N \parSpace_n, \quad \parSpace_n =[a_n,b_n],
\]
for some $N\in\N$.
Extensions to cases in which $\parSpace_n$ are endowed with a probability measure (as typical in uncertainty quantification) and possibly unbounded intervals are straightforward.
Furthermore, let $q: \solSpace \to \R^T$ be a map from the solution space of \cref{eq:PDE} to $\R^T$,
chosen to extract quantities of relevance from the PDE solution $\sol$,
such as the flux across a prescribed surface in fluid dynamics,
stress and strain fields in structural mechanics, and drag or lift forces in aerodynamics.
The \emph{Quantity of Interest} (\QoI) is then defined as the composition 
$\qoi:\parSpace\to\R^T$,
\begin{equation}\label{eq:qoidef}
  \qoi(\parVar) \coloneqq \qoisol \left(\sol\right),  
\end{equation}
which encodes the dependence of the output on the parameters $\parVar$.
In practice, only a numerical approximation $\asol$ of $\sol$ is available,
\eg obtained by finite elements;
for the rest of this work, with a slight abuse of notation we use $\qoi(\parVar)$ to denote also $\qoisol \left(\asol\right)$.
Furthermore, in the following we consider for simplicity scalar-valued \QoI ($T=1$); extensions to vector-valued \QoI are straightforward.

As already mentioned, we are interested in analyzing how $\qoi$ changes as a function of the parameters $\parVar$,
which is typically unaffordable by brute-force inspection, since evaluating $\qoi(\parVar)$
requires computing $\sol$
and is thus an expensive operation.
To mitigate this problem, we build a \emph{surrogate model} of $\qoi$, \ie an approximation of $\qoi(\parVar)$ (\eg a polynomial)
employing a relatively small set of values $\qoi(\parVar_1),\ldots,\qoi(\parVar_M)$. 
The resulting surrogate model is much cheaper to evaluate, 
and it can be queried multiple times for a limited computational cost.

Obtaining values of $\qoi(\parVar)$ means in particular solving the non-linearity in \cref{eq:PDE},
which is typically done by iterative methods that 
require the solution of a linear system at each
step, \eg the Newton's method
\begin{equation}\label{eq:newton-iter}
\sol^{\{i+1\}} = \sol^{\{i\}} - \big(\delta_{\solSpace}\diffOp_{(\parVar, \sol^{\{i\}})}\big)^{-1} \diffOp(\parVar, \sol^{\{i\}})
\end{equation}
where $\delta_{\solSpace}\diffOp_{(\parVar, \sol^{\{i\}})}$
is a linear operator, namely the Fréchet derivative
of $\diffOp$ with respect to its functional argument
evaluated at $(\parVar, \sol^{\{i\}})$.

Formulas for the partial derivatives of $\sol$ and $\qoi$ with respect to the parameters are instead obtained starting from the  derivative of \cref{eq:PDE} with respect of $\uparVar_n$, that is $0$ since \cref{eq:PDE} holds for all $\parVar$, giving
\[
  \partial_{\uparVar_n} \diffOp (\parVar, \sol ) +
  \delta_{\solSpace}\diffOp_{(\parVar,\sol)} \partial_{\uparVar_n}\sol = 0,
\]
where $\partial_{\uparVar_n} \diffOp$ is the partial derivative of $\diffOp$ with respect of the $n$-component of its first argument. From the previous equation,
one gets
\begin{align}
  \partial_{\uparVar_n}\sol & = -\big(\delta_{\solSpace}\diffOp_{(\parVar,\sol)}\big)^{-1} \partial_{\uparVar_n}\diffOp(\parVar,\sol) \label{eq:par-derivative-of-u}, \\
  \partial_{\uparVar_n}\qoi(\parVar) & = 
  \delta_{\solSpace}\qoisol_{\sol} \partial_{\uparVar_n}\sol 
  = -\delta_{\solSpace}\qoisol_{\sol} \big(\delta_{\solSpace}\diffOp_{(\parVar,\sol)}\big)^{-1}\partial_{\uparVar_n} \diffOp (\parVar, \sol), \label{eq:par-derivative-of-qoi}
\end{align}
where the expression for $\partial_{\uparVar_n}\qoi(\parVar)$ is obtained by applying the chain rule of derivation to
\cref{eq:qoidef} and $ \delta_{\solSpace}\qoisol_{\sol}$
is the Fréchet derivative of $\qoi$ with respect to its functional argument evaluated at $\sol$.
From the equations above we deduce that, once Newton's method has converged,
computing $\partial_{\uparVar_n}\qoi(\parVar)$ costs exactly
one linear solve and is thus not more expensive than computing $\qoi(\parVar)$ for a new $\parVar \in \parSpace$
(the two operations have the same cost whenever the Newton's method converges in one iteration).
The same conclusion is still valid if one uses finite differences instead of \cref{eq:par-derivative-of-qoi}
to approximate $\partial_{\uparVar_n}\qoi$, since such computation would require two evaluations of $\sol$, one of which is already available.
For generic solvers one can resort to \emph{algorithmic differentiation} \cite{Griewank2003} leading to similar conclusions.
Thus, depending on the number of parameters and on the number of Newton iterations,
it can be computationally cheaper to compute $\nabla_{\parSpace}\qoi$ at a $\parVar\in\parSpace$
where $\sol$ has already been computed rather than computing $\qoi$ at a different $\parVar \in \parSpace$,
which makes it attractive to devise a mechanism to use $\nabla_{\parSpace}\qoi$
in the construction of a surrogate model. 

\begin{remark}\label{rem:assembly-costs}
  The difference in computational cost between computing $\nabla_{\parSpace}\qoi$ and
  evaluating $\qoi(\parVar)$ at a new point $\parVar \in \parSpace$ go potentially
  beyond the mere count of Newton iterations.
  To begin with, computing the full gradient $\nabla_\parSpace\qoi$
  amounts to solving \cref{eq:par-derivative-of-u} for each $n = 1, \dots, N$,
  \ie solving $N$ linear systems with the \emph{same} operator
  $\delta_{\solSpace}\diffOp_{(\parVar,\sol)}$ for all $n$
  (which can therefore be factorized once and reused)
  and $N$ different right-hand sides $\partial_{\uparVar_n}\diffOp(\parVar,\sol)$.
  By contrast, each Newton iteration \cref{eq:newton-iter} requires \emph{one} linear solve,
  but with a matrix that must be \emph{updated} at every step.
  Furthermore, for a scalar-valued \QoI or a vector-valued one with $T<N$,
  by taking a closer look at the vector version of \cref{eq:par-derivative-of-qoi}, \ie
\[
  \nabla_{\parSpace}\qoi(\parVar)
  = -\delta_{\solSpace}\qoisol_{\sol} \big(\delta_{\solSpace}\diffOp_{(\parVar,\sol)}\big)^{-1}\nabla_{\parSpace} \diffOp (\parVar, \sol),
\]
one can see that the cost of computing $\nabla_{\parSpace}\qoi$ can be further
reduced to that of $\min\{T,N\}=T$ linear solves, by computing as a first step the
product $\delta_{\solSpace}\qoisol_{\sol}
\big(\delta_{\solSpace}\diffOp_{(\parVar,\sol)}\big)^{-1}$, \ie by solving the adjoint problem of the generic Newton iteration.
\end{remark}

\section{Surrogate models via polynomial interpolation using derivative information}\label{sec:hermite}

Polynomial interpolation represents one of the simplest 
methods of surrogate modeling.
In the following, we recall the main aspects of derivatives-informed polynomial interpolation in both the univariate and multivariate settings.
Specifically, we focus on the case in which at each point of $\parSpace$
the values of function and derivatives up to a certain prescribed order
(possibly different at each point) are required to match. 
\begin{remark}
  In this work, we use the term \emph{Hermite interpolation} to refer to interpolation problems involving function values and possibly higher-order derivatives.
  We use instead the term \emph{osculatory interpolation} when only function values and first order derivatives are involved.
  Part of the literature uses the opposite convention, and there is no general fixed terminology for this topic. 
\end{remark}

\subsection{Univariate interpolation}\label{sec:univariate-interpolation}
If $\parSpace=[a,b]\subset\R$, the general Hermite interpolation problem reads:
given a function $\qoi:\parSpace\to\R$, a set of $M$ interpolation nodes
$\interpPoints=\{\interpPoint_1,\dots,\interpPoint_M\}\subset \parSpace$,
a set of $M$ non-negative integers $\bar k_1,\dots, \bar k_M$, and a set of observations
\[
  \partial^k \qoi(\interpPoint_j),\qquad \forall j=1,\dots,M,\ k=0,\dots, \bar k_j,
\]
find
$\aqoi\in \poly$ such that
\[
  \partial^k\aqoi(\interpPoint_j)=\partial^k\qoi(\interpPoint_j),\qquad \forall j=1,\dots,M,\ k=0,\dots, \bar k_j;
\]
where $\bar k_1,\dots, \bar k_M$ specify the maximum derivative order prescribed at each node. 
Therefore, the total number of conditions imposed is $M+\sum_{j=1}^M \bar k_j$.
Upon choosing $\poly$ as the space of polynomials of degree $\le \deg$, 
the problem is well-posed if and only if 
\begin{equation}\label{eq:dof1D}
  \dim \poly = \deg + 1 = M + \sum_{j=1}^M \bar k_j.
\end{equation}
For a well-posed problem, let $\lagrange_{j,k}$ be the unique element of $\poly$ satisfying
\begin{equation}\label{eq:double-delta}
\partial^s\lagrange_{j,k} (\interpPoint_i)=\delta_{j,i}\delta_{s,k},\qquad \forall i,j=1,\dots,M,\ s,k=0,\dots, \bar k_j.  
\end{equation}
The set $\{\lagrange_{j,k}\}\subset\poly$ is the \emph{nodal basis} of $\poly$ associated with the interpolation problem,
and the natural generalization of the classical Lagrange basis to the Hermite setting.
The solution of the interpolation problem can be computed explicitly and is given by the operator $\U$ defined as follow
\[
\U\qoi(\uparVar)\coloneqq\sum_{j=1}^M \sum_{k=0}^{\bar k_j} \partial^k\qoi(\interpPoint_j) \lagrange_{j,k}(\uparVar).
\]
In the following, we provide additional details for Lagrange (values only) and osculatory (values and first order derivatives) polynomial interpolation.

\subsubsection*{Lagrange interpolation}
The case $0=\bar k_1=\dots=\bar k_M$ is known as Lagrange
interpolation, and is well-posed for $\deg=M-1$. In this setting, the nodal basis functions $\lagrange_{j,0}$
are the standard Lagrange polynomials $\lagrange_j$ (see \cref{fig:oscpol}-left)
and we denote the interpolation operator by $\L$, defined as follow
\begin{equation}\label{eq:unilag}
\L\qoi(\uparVar)\coloneqq \sum_{j=1}^M \qoi(\interpPoint_j) \lagrange_{j}(\uparVar),\quad \lagrange_{j}(\uparVar)\coloneqq\prod_{i\ne j}\frac { \uparVar-\interpPoint_i}{ \interpPoint_j-\interpPoint_i}.
\end{equation}

\subsubsection*{Osculatory interpolation} 
\begin{figure}[!t]
  \centering
  \hspace{-1em}\includegraphics{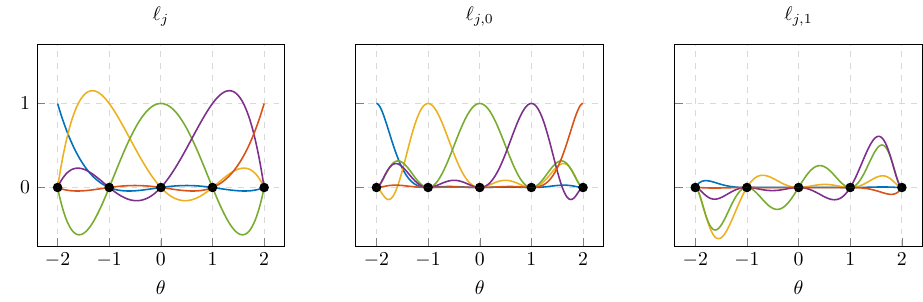}
  \caption{Lagrange polynomials $\lagrange_j$ (left), compared with osculatory polynomials $\lagrange_{j,0}$ (center) and $\lagrange_{j,1}$ (right),
    for $j=1, \ldots, M$, associated to the nodes $\{-2, -1, 0, 1, 2\}$.}
  \label{fig:oscpol}
\end{figure}

The case $1=\bar k_1=\dots=\bar k_M$ is well-posed for $\deg=2M-1$.
The interpolation operator in this setting is defined as 
\begin{align}\label{eq:uni-interp-osc}
&\O\qoi(\uparVar)\coloneqq\sum_{j=1}^M \Big(  \qoi(\interpPoint_j) \lagrange_{j,0}(\uparVar)+ \partial\qoi(\interpPoint_j) \lagrange_{j,1}(\uparVar)\Big),
\\\notag&\text{where}\quad \lagrange_{j,0}(\uparVar)= \Big(1-2(\uparVar-\interpPoint_j)  \lagrange_j'(\interpPoint_j )\Big)\lagrange_j(\uparVar )^2,\,\qquad
\lagrange_{j,1}(\uparVar)= (\uparVar-\interpPoint_j) \lagrange_j(\uparVar )^2.
\end{align}
As stated in \cref{eq:double-delta}, 
each $\lagrange_{j,0}$ has values $1$ on interpolation nodes $\interpPoint_j$ and
$0$ on the other nodes, and its first derivative is $0$ on all interpolation nodes
(see \cref{fig:oscpol}-center).
Conversely, each $\lagrange_{j,1}$ has value $0$ on all interpolation nodes,
and its first derivative has values $1$ on interpolation nodes $\interpPoint_j$
and zeros on the other nodes (see \cref{fig:oscpol}-right).

\subsubsection*{Comparing Lagrange and osculatory interpolation}
A fair way of comparing the approximation error of a function $\qoi$
via Lagrange and osculatory interpolation is by ensuring that the two
approximants belong to the same space, \ie they have the same degree $\deg$.
Since each interpolation point corresponds to
two data for osculatory interpolation, the comparison should
thus necessarily use two sets of points with different cardinality,
\ie $\interpPoints_{\L}$ with $M_{\L}=\deg+1$ points for Lagrange interpolation
and  $\interpPoints_{\O}$ with $M_{\O}=M_{\L}/2=(\deg+1)/2$ for osculatory interpolation,
with $\deg$ odd.

It is well-understood that the approximation properties
of any interpolant depend on the placement of
the interpolation points: badly-positioned points make the interpolation
operator unstable. This behavior can be analyzed in term of the Lebesgue constant
of the interpolation problem, \ie the norm of the interpolation operator.
Lebesgue constants have been studied for both Lagrange and
osculatory interpolation \cite{manni1993,Szabados1993357}.
Comparing Lagrange and osculatory interpolants by comparing their
Lebesgue constants is however misleading because
they measure norms of operators between different spaces: indeed, the domain of the Lagrange interpolant operator is the space of
  \emph{continuous} functions $C^0\left(\parSpace\right)$, whereas the domain of
  the osculatory interpolant is the space of function with
  \emph{continuous derivatives} $C^1\left(\parSpace\right)$,
  \ie a smaller class of functions endowed with a stronger norm.
Therefore, we must follow a different approach.

While a complete discussion about the comparison between the two interpolation schemes
is beyond the scope of this work, sharp a-priori error estimates for regular functions can be obtained from the residual formula
of interpolation and allow to compare the two approaches on the same class of functions, namely $\qoi$
such that $\qoi^{(\deg+1)}$ is bounded. The residual formula for Lagrange polynomial
interpolation is: for all $\uparVar\in\parSpace$ there exists
$\zeta\in\parSpace$ such that
\begin{equation}\label{eq:interpolation-error}
  (\qoi-\L\qoi )(\uparVar)=  \frac{\qoi^{(\deg+1)}(\zeta)}{(\deg+1)!}
  \prod_{j=1}^{M_{\L}} (\uparVar -\interpPoint_j),
\end{equation}
from which we get directly the error estimate for Lagrange interpolation
\begin{equation*}\label{eq:lagrange-interpolation-error}
	\norm{\qoi-\L\qoi}\le  \frac{C_{\interpPoints_{\L}}}{(\deg+1)!}  \norm{\qoi^{(\deg+1)}},\qquad C_{\interpPoints_{\L}}=\norm{\prod_{j=1}^{M_{\L}} (\uparVar -\interpPoint_j)}.
\end{equation*}
The constant $C_{\interpPoints_{\L}}$ is minimized when
$\interpPoints_{\L} = \interpPoints_{\L}^\star$
is the set of zeros of the Chebyshev polynomial of degree $\deg+1$
(see \eg \cite[Chapter 2, Theorem 11]{lorentz:book}), in which case we have
\begin{equation}\label{eq:lagrange-interpolation-error-best-points}
C_{\interpPoints_{\L}^\star}=2^{-2d-1}|\parSpace|^{\deg+1}.
\end{equation}
Next, we exploit the fact that osculatory interpolation 
on a set of points, $\interpPoints_{\O}$, can be seen as the limit of 
Lagrange interpolation on the set $\interpPoints_{\L}$ defined as follows: 
$\interpPoints_{\L} = \interpPoints_{\O} \cup (\interpPoints_{\O}+\epsilon)$
when $\epsilon\to 0$, \ie when each pair of $\epsilon$-separated points
collapses to a single point:
$\interpPoint_{2j}=\interpPoint_{2j-1}$ for $j=1,\dots,M_{\O}$.
The reminder formula \cref{eq:interpolation-error} applies giving
the following specialized version
\[
  (\qoi-\O\qoi )(\uparVar)=  \frac{\qoi^{(\deg+1)}(\zeta)}{(\deg+1)!}
  \prod_{j=1}^{M_{\O}} (\uparVar -\interpPoint_{2j})^2,
\]
from which it follows
\begin{equation*}\label{eq:osculatory-interpolation-error}
  \norm{\qoi-\O\qoi}\le  \frac{C_{\interpPoints_{\O}}^2}{(\deg+1)!}  \norm{\qoi^{(\deg+1)}},
  \qquad C_{\interpPoints_{\O}}=\norm{\prod_{j=1}^{M_{\O}} (\uparVar -\interpPoint_j)}.
\end{equation*}
The above is minimized choosing $\interpPoints_{\O}=\interpPoints_{\O}^\star$
as the zeros of the Chebyshev polynomial of degree $(\deg+3)/2$
in which case we have
\begin{equation}\label{eq:osculatory-interpolation-error-best-points}
C_{\interpPoints_{\O}^\star}^2=2^{-2d}|\parSpace|^{\deg+1}.
\end{equation}
Remembering that the constant in \cref{eq:lagrange-interpolation-error-best-points} 
and \cref{eq:osculatory-interpolation-error-best-points} are sharp, \ie that
there exists a function for which the estimates hold with equality
(namely, the monomial $\uparVar\mapsto\uparVar^{\deg+1}$), we conclude that
with optimal point-placement the a-priori error estimate for osculatory interpolation is
twice as large as that of Lagrange interpolation \emph{for the same polynomial space}.
However, if the cost of obtaining an evaluation of $\partial \qoi$
is sufficiently smaller than the cost of obtaining an evaluation of $\qoi$,
as discussed in \cref{sec:parametricPDEs},
it is possible that an osculatory interpolant provides a better error
than Lagrange interpolation \emph{for the same computational cost},
also for suboptimal choices of point-placement.

\subsection{Tensor-product interpolation}\label{sec:multivariate-interpolation}
Tensor-product constructions are a straightforward way to extend univariate approximations to the multivariate setting,
allowing the construction of well-posed schemes.
Specifically, tensor-product interpolation reduces multivariate interpolation problem to
many univariate interpolations by imposing additional structure on the set of
interpolation points $\interpPoints$ and on the function space $\poly$.

Recalling $\parSpace=\bigotimes_{n=1}^N \parSpace_n$, the idea is to consider a univariate interpolation scheme per dimension: for each $n=1, \ldots, N$, let
\begin{align*}
  &\interpPoints^{[n]}=\{\interpPoint^{[n]}_{1}, \ldots,\interpPoint^{[n]}_{M_n}\}\subseteq \parSpace_n,\\
  &\bar k^{[n]}_{1}, \ldots, \bar k^{[n]}_{M_n} \in\N, \\
  &\poly^{[n]} \text{ a space of polynomials } \parSpace_n\to\R.
\end{align*} 
Then, given $\qoi: \parSpace\to\mathbb{R}$, the interpolation problem on the Cartesian grid
$\interpPoints=\interpPoints^{[1]}\times\dots\times\interpPoints^{[N]}$,
with the tensor-product space $\poly=\poly^{[1]}\otimes\dots\otimes\poly^{[N]}$ is: find $\aqoi\in \poly$
such that 
\begin{equation}\label{eq:tp-hermite-interpolation}
  \partial^{k_1}_{\uparVar_1}\ldots\partial^{k_N}_{\uparVar_N}\aqoi(\bm \interpPoint_{\bm j})
  =
  \partial^{k_1}_{\uparVar_1}\ldots\partial^{k_N}_{\uparVar_N}\qoi(\bm \interpPoint_{\bm j}),
  \quad \forall \bm \interpPoint_{\bm j}=(\interpPoint^{[1]}_{j_1},\dots,\interpPoint^{[N]}_{j_N})\in \interpPoints,\  0\le k_n\le \bar k^{[n]}_{j_n}. 
\end{equation}

The problem is well-posed if and only if for each $n=1, \ldots, N$ the
corresponding univariate Hermite interpolation problem is well-posed.
The solution of \cref{eq:tp-hermite-interpolation} can be expressed as
\begin{equation*}\label{eq:tensor-product-interpolation}
	\U\qoi(\parVar)=\sum_{\bm \interpPoint_{\bm j} \in \bm\interpPoints} \sum_{\bm 0\le \bm k\le  \bar{\bm k}_{\bm \interpPoint}}\partial^{\bm k} \qoi(\bm \interpPoint_{\bm j})\  \lagrange_{\bm j, \bm k}(\bm \uparVar),\qquad \lagrange_{\bm j, \bm k}(\bm \uparVar) = \prod_{n=1}^{N}\lagrange^{[n]}_{j_n, k_n}(\uparVar_n),
\end{equation*}
where $\lagrange^{[n]}_{j_n, k_n}$ are the nodal basis functions described in the previous subsection.
Tensor-product constructions are convenient for their simplicity,
but they are practical only for very small number of dimensions $N$,
due to the well-known phenomenon of \emph{curse of dimensionality} \cite{bellman1961}.

\subsubsection*{Tensor-product Lagrange interpolation}
The case $\bar k^{[n]}_j=0$ for $n=1,\dots,N$ and $j=1,\dots,M_n$, \ie
  the case in which the
values of $\qoi$ are interpolated on a Cartesian grid, manifests the
\emph{sampling} curse of dimensionality: for a grid with $M_n = \bar{M}$
points per direction, the number of interpolation points is indeed $\bar{M}^N$
\ie it increases exponentially with $N$.
The Smolyak algorithm described in \cref{sec:abstractsmolyak} is
a classical strategy to circumvent this issue.

\subsubsection*{Tensor-product osculatory interpolation} 
A closer look to \cref{eq:tp-hermite-interpolation} shows that
the case $\bar k^{[n]}_j=1$ for $n=1,\dots,N$ and $j=1,\dots,M_n$ interpolates not only the partial derivatives $\partial_{\uparVar_n}\qoi$
at the interpolation nodes, but also all the mixed derivatives of order up to $1$ in all directions,
thus requiring the evaluation of $2^N$ derivatives of $\qoi$ at each interpolation point.
For example, constructing a bivariate osculatory interpolant on a $5\times 5$ grid (\ie $N=2$, $M_1=M_2=5$, $\bar k^{[1]} = \bar k^{[2]} =1$)
requires determining $10 \times 10 = 100$ degrees of freedom, cf. equation \cref{eq:dof1D}.
The evaluations of $\qoi$ at the $25$ grid points provide $25$ data, and the first order derivatives $\partial_{\uparVar_1} \qoi$ and $\partial_{\uparVar_2} \qoi$
give us additional $50$ data, for a total of $75$ data; the remaining degrees of freedom are fulfilled upon knowing the 
$25$ mixed derivatives $\partial_{\uparVar_1}\partial_{\uparVar_2}\qoi$.

This interpolation hence suffers twice from the curse of dimensionality: once due the number of grid points, \ie the sampling curse of dimensionality discussed earlier,
  and once due to the amount of information that needs to be computed at each grid point,
  that we name \emph{derivative-induced} curse of dimensionality.
Moreover, in practical scenarios where $\qoi$ is a \QoI coming from a PDE, any information beyond first derivatives $\partial_{\uparVar_n} \qoi$ is essentially never available, 
making this interpolation scheme unsuitable for applications.

\section{A sparse-grids-inspired approach to multi-variate osculatory approximation}\label{sec:SG}

A common tool in literature to address the two above-mentioned curses of dimensionality is the Smolyak algorithm, recalled in abstract terms
in \cref{sec:abstractsmolyak} following \cite{bungartz2004sparse,wolfers.tempone:smolyak}. 
In the context of Lagrange interpolation, it leads to replacing Cartesian grids with
sparse grids \cite{babuska.nobile.eal:stochastic2, bungartz2004sparse, xiu.hesthaven:high},
see \cref{sec:sparse-grids}.
This algorithm has however a broader applicability and can be invoked
whenever one has a computational problem with a tensor-product structure
whose cost needs to be reduced: in particular, we show in
\cref{sec:osculatory-smolyak} that in the context of osculatory interpolation
it constructs a osculatory interpolant that employs only gradient information
($N$ unidirectional partial derivatives) instead of all of the
$2^N$ mixed partial derivatives on a Cartesian grid.
Unfortunately, the two approaches cannot be used together,
as shown in \cref{sec:challenges}:
in \cref{sec:GELS} we therefore propose a different remedy
to the two curses of dimensionality, combining a sparse-grids sampling scheme and a least-squares approach to approximation.

\subsection{Abstract Smolyak algorithm}\label{sec:abstractsmolyak}

For each direction $\uparVar_n$, $n=1, \ldots, N$
we consider a sequence of approximation operators $0=\U_{0}^{[n]}$, $\U_{1}^{[n]}$, $\dots$, $\U_{i_n}^{[n]}$, $\dots$ of increasing accuracy.
We refer to $i_n\in \N$ as \emph{univariate level} and we denote the range of $\U_{i_n}^{[n]}$ with $\polyop{\U_{i_n}^{[n]}}$.
To fix ideas, we already mention that
in the following either $\U_{i_n}^{[n]}$ are Lagrange interpolation operators whose number of interpolation points increases with $i_n$ (see \cref{sec:sparse-grids}), or
$\U_{1}^{[n]}$ is a Lagrange interpolant and $\U_2^{[n]}$ is an osculatory interpolant
on the same number of interpolation points (see \cref{sec:osculatory-smolyak}).

To move to the multi-variate setting, let $\bm i = (i_1, \ldots, i_N) \in \N^N$
be a multi-index that collects univariate levels along each direction. We then define the tensor-product operator
\[
  \U_{\bm i} \coloneqq \U_{i_1}^{[1]} \otimes \dots \otimes \U_{i_N}^{[N]},
\]
whose range is the tensor-product space
\[
  \polyop{\U_{\bm i}} \coloneqq \polyop{\U_{i_1}^{[1]}} \otimes \dots \otimes \polyop{\U_{i_N}^{[N]}}.
\]
Therefore $\U_{\bm i}\qoi \in \polyop{\U_{\bm i}}$ is  an approximation of $\qoi$ that employs $\prod_{n=1}^N \dim \polyop{\U_{i_n}^{[n]}}$
degrees of freedom (where the curse of dimensionality shows up)
and gets increasingly accurate if $i_n \rightarrow \infty$ in every direction $n$. 

The Smolyak algorithm proceeds by defining the univariate detail operators for $i_n \in \Np$
\[
  \Delta_{i_n}^{[n]}\coloneqq \U_{i_n}^{[n]}-\U_{i_n-1}^{[n]},
\]
and multivariate detail operators for $\bm i \in \Np^N$
\begin{equation}\label{eq:multivar-delta-def}
  \Delta_{\bm i}\coloneqq \Delta_{i_1}^{[1]}\otimes \dots \otimes \Delta_{i_N}^{[N]}
  = \sum_{\bm j \in \{0,1\}^N}(-1)^{\|\bm j\|_1}\U_{\bm i - \bm j}.  
\end{equation}
Note in particular that the following telescopic equality holds:
\begin{equation}\label{eq:telescopic}
  \U_{\bm r} = \sum_{ \bm i \in \indexset_{\text{TP},\bm r}}\Delta_{\bm i}, \quad
  \indexset_{\text{TP},\bm r} = \{ \bm i \in \Np^N :  i_n \leq r_n \},
\end{equation}
showing that the tensor-product operator $\U_{\bm r}$ can be decomposed  as a sum of multivariate detail operators.
For any multi-index set $\indexset\subset\Np^N$,
we finally define the Smolyak operator 
\begin{equation}\label{eq:sg-operator}
  \S_\indexset \coloneqq \sum_{\bm i \in \indexset}\Delta_{\bm i}.
\end{equation}
The underlying idea is that generalizing \cref{eq:telescopic} to \cref{eq:sg-operator} by changing $\indexset_{\text{TP},\bm r}$
  with $\indexset$
  one can discard the terms in \cref{eq:telescopic} whose contribution to the overall approximation
$\U_{\bm r}\qoi$ is negligible; we  return on the choice of $\indexset$ later on.
Using \cref{eq:multivar-delta-def}, $\S_\indexset$ can be further re-expressed as a linear combination of tensor-product operators (with integer coefficients), namely
\begin{equation}\label{eq:ct}
  \S_\indexset = \sum_{\bm i \in \overline{\indexset}}c_{\bm i}\U_{\bm i}, \quad c_{\bm i} = \sum_{\substack{\bm j\in\{0,1\}^N,\\ \bm i + \bm j \in\overline{\indexset}}}(-1)^{\|\bm j\|_1},
\end{equation}
where $\overline\indexset\coloneqq\{\bm i-\bm j:\bm i \in \indexset,\, \bm j \in \{0,1\}^N\}$;
note that $c_{\bm i}$ can be zero, in which case the corresponding operators
$\U_{\bm i}$ do not need to be assembled in practice. 
We refer to \cref{eq:ct} as the \emph{combination-technique} form of the Smolyak operator,
and to the tensor-product operators $\U_{\bm i}$ that enter the linear combination
with $c_{\bm i} \neq 0$ as \emph{component operators}.
Given \cref{eq:ct} it is easy to see that the Smolyak operator delivers an approximation in the composite space
\begin{equation}\label{eq:polysparse}
\polyop{\S_\indexset} \coloneqq \sum_{\bm i \in\overline{\indexset},\,c_{\bm i} \neq 0}\polyop{\U_{\bm i}}.  
\end{equation}

As already mentioned, the efficiency of the Smolyak operator,
\ie its ability to deliver an approximation $\S_\indexset\qoi$
that is close to its tensor-product operator counterpart $\U_{\bm r}\qoi$
for a fraction of the degrees of freedom (hence computational cost),
crucially depends on the choice of the multi-index set $\indexset$.
A classical choice to this end is 
\begin{equation}\label{eq:TDset}
  \indexsetsum=\left\{\bm i \in \Np^N: \sum_{n=1}^N (i_n-1) \le w \right\},
\end{equation}
for some $w \in \N$ that in this work is called \emph{Smolyak level};
\cref{eq:TDset} is actually the multi-index set that we employ in this work;
note that $\indexsetsum$ is such that $\indexset = \overline{\indexset}$.
We refer to \cite{bungartz2004sparse,wolfers.tempone:smolyak} for a wider discussion on the design of $\indexset$ either by
a-priori strategies or by adaptive a-posteriori algorithms inspired by the seminal work \cite{gerstner.griebel:adaptive}.
In particular, if $\qoi$ is regular enough then $\S_{\indexsetsum}\qoi \rightarrow \qoi$ in a suitable sense 
as $w \rightarrow \infty$, see again \eg \cite{bungartz2004sparse,wolfers.tempone:smolyak}.

\subsection{Smolyak algorithm for the sampling curse of dimensionality}\label{sec:sparse-grids}

For each parameter $\uparVar_n$ and univariate level $i_n\in \N$, consider the Lagrange interpolation operator
$\L_{i_n}^{[n]}$ defined in \cref{eq:unilag} on a set
$\interpPoints_{i_n}^{[n]}\subset\parSpace_n$ of $m(i_n)$ points, where
$m:\N\to\N$ with $m(0)=0$ is a non-decreasing function that prescribes the
number of interpolation points for the each univariate level. 
Next, let $\L_{\bm i}\coloneqq \L^{[1]}_{i_1}\otimes\dots\otimes\L^{[N]}_{i_N}$ be the Lagrange interpolant on the Cartesian grid
$\interpPoints_{\bm i}\coloneqq \interpPoints_{i_1}^{[1]}\times \ldots \times \interpPoints_{i_N}^{[N]}$ that contains $\prod_{n=1}^N m(i_n)$ points, showing the \emph{sampling} curse of dimensionality.

We denote by $\SL_{\indexset}$ the Smolyak operator $\S_\indexset$ defined in
\cref{eq:ct} when $\U_{\bm i}=\L_{\bm i}$, and given our choice $\indexset=\indexsetsum$ as in \cref{eq:TDset} we introduce the compact notation
\begin{equation}\label{eq:lagr-sp-grid}
\SL_{w}\coloneqq \SL_{\indexsetsum} = \sum_{\bm i \in \indexsetsum}c_{\bm i}\L_{\bm i}.  
\end{equation}
In general, for any multi-index set $\indexset$ we refer to the grids $\interpPoints_{\bm i}$ with $\bm i\in \overline{\indexset}$ and $c_{\bm i} \neq 0$ in \cref{eq:ct}
  as the \emph{component grids} and to their union as \emph{sparse grid}: 
\begin{equation*}\label{eq:sg-cartesian-grid}
  \interpPoints_{\indexset} \coloneqq
  \bigcup_{\bm i \in\overline{\indexset},\,c_{\bm i} \neq 0}\interpPoints_{\bm i}.
\end{equation*}
Note that the approximation
$\SL_{\indexset}\qoi$ is constructed using the values of $\qoi$ on
$\interpPoints_{\indexset}$: the location of the interpolation
points $\interpPoints_{i_n}^{[n]}$ plays a fundamental role for both the quality of the final
approximation $\SL_{\indexset}\qoi$ and its implementation efficiency.
In particular \emph{nested} interpolation points,
$\interpPoints_{j_n}^{[n]} \subset \interpPoints_{i_n}^{[n]}$ if $j_n<i_n$, allow reusing the evaluations of $\qoi$
for the construction of multiple component operators of $\SL_{\indexset}$
and is therefore computationally advantageous.
Moreover, in this context, it can be proved that $\SL_{\indexset}$ 
is an interpolant on $\interpPoints_{\indexset}$
if and only if the univariate points $\interpPoints_{j_n}^{[n]}$ are nested.
For this reason, $\SL_{\indexset}$ is typically called \emph{sparse-grid interpolant}.

\begin{figure}[!t]
	\centering
	\includegraphics{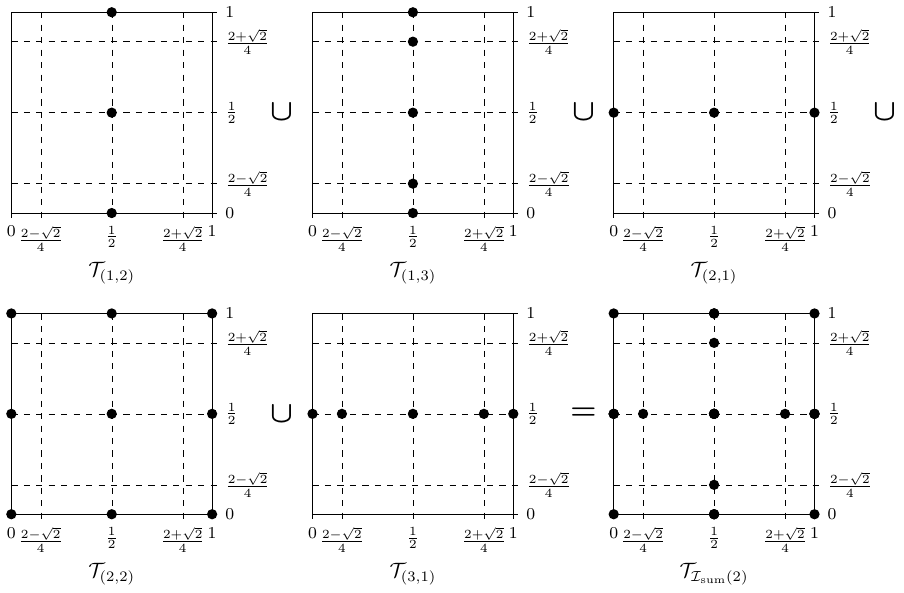}
	\caption{Sparse grid in $N=2$. $\U_{i_n}^{[n]}$ employs Clenshaw--Curtis points \cref{eq:CCpoints}, doubling-rule \cref{eq:m-doubling} for $m(i_n)$,
and $\indexset = \indexsetsumw{2}  =  \{(1,2), (1,3), (2,1), (2,2), (3,1)\}$, cf. \cref{eq:TDset}.
          The panels show the different component grids $\mathcal{T}_{\bm i}$ and the final sparse grid $\mathcal{T}_{\indexsetsumw{2}}$.}
	\label{fig:sg-plot}
\end{figure}

The specific form of the polynomial space $\polyop{\SL_{\indexset}}$ to which a sparse-grid approximation $\SL_{\indexset}\qoi$ belongs, see  \cref{eq:polysparse},  depends on the choice of
the multi-index set $\indexset$ and on the function $m$ see \cite{back2010} for a general discussion. 
In this work we set $m$ to be the \emph{doubling rule}
\begin{equation}
  m(i_n) =
   \begin{cases}
     i_n & i_n=1,\\
     2^{i_n-1}+1 & i_n>1,
   \end{cases} \label{eq:m-doubling}
\end{equation}
and together with \cref{eq:TDset} we get (see \cite{back2010})
\begin{equation}\label{eq:sparse-grid-poly-space}
  \begin{aligned}
  \polyop{\L_{\bm i}} & = \spn\left\{\uparVar_1^{j_1}\dots\uparVar_N^{j_N} : j_n\in \N,\, 0 \leq j_n \leq m(i_n)-1  \right\}, \\   
  \polyop{\SL_{w}}
  & = \spn\left\{\uparVar_1^{j_1}\dots\uparVar_N^{j_N} : j_n\in \N,\, \sum_{n=1}^N g(j_n) \le w\right\}, \quad
g(j_n)  =
  \begin{cases}
    j_n & j_n=0,1, \\
    \lceil \log_2(j_n) \rceil & j_n \geq 2.
  \end{cases}
\end{aligned}
\end{equation}
\cref{fig:sg-plot} shows the sparse grid $\mathcal{T}_{\indexsetsumw{2}}$ and its component grids $\mathcal{T}_{\bm i}$, for $N=2$ and Clenshaw--Curtis points \cite{trefethen2008}: 
\begin{equation} \label{eq:CCpoints}
  \interpPoint_\kappa= \cos\left( \frac{(\kappa-1)\pi}{m(i_n)-1}\right), \quad 1 \leq \kappa \leq m(i_n).
\end{equation}
Using the doubling rule \cref{eq:m-doubling}  with Clenshaw--Curtis points is a classical choice that guarantees nestedness of the interpolation points $\interpPoints_{i_n}^{[n]}$ and is employed throughout this work; more choices are discussed, \eg in \cite{piazzola2024}.

\subsection{Smolyak algorithm for the derivative-induced curse of
  dimensionality}\label{sec:osculatory-smolyak} 
Within the setting of the previous section,
consider now a fixed grid $\interpPoints_{\bm r}$ with $\prod_{n=1}^N m(r_n)$ points.
Let
\[
  \U_{1}^{[n]}=\L_{r_n}^{[n]},\qquad \U_{2}^{[n]}=\O_{r_n}^{[n]},
\]
where $\L_{r_n}^{[n]}$ is the Lagrange interpolant on $\interpPoints_{r_n}^{[n]}$,
and $\O_{r_n}^{[n]}$ is the osculatory interpolant, see \cref{eq:uni-interp-osc},
on the same set of points.
Since we are not interested in interpolating higher-order derivatives,  
we do not define $\U_{i_n}^{[n]}$ for $i_n>2$,
and consequently we apply the Smolyak algorithm \cref{eq:ct} with the
multi-index set $\indexsetsumw{1} = \{ \bm 1, \dd_1, \ldots, \dd_N\}$, where $\bm e_n$ is
the multi-index whose $n$-th component is 1 and the remaining others are 0. 
We thus obtain the component operators
\[
\U_{\bm 1}=\L_{\bm r},\quad\text{and}\quad \U_{\dd_n}=\L_{r_1}^{[1]}\otimes\dots\otimes\L_{r_{n-1}}^{[n-1]} \otimes\O_{r_{n}}^{[n]}\otimes\L_{r_{n+1}}^{[n+1]}\otimes\dots\otimes\L_{r_N}^{[N]},
\]
and we denote the resulting Smolyak operator as
\begin{equation}\label{eq:ct-single-grid}
  \SO_{\bm r} \coloneqq
  - (N-1) \L_{\bm r} + \sum_{n=1}^N \L_{r_1}^{[1]}\otimes\dots\otimes\L_{r_{n-1}}^{[n-1]} \otimes\O_{r_{n}}^{[n]}\otimes\L_{r_{n+1}}^{[n+1]}\otimes\dots\otimes\L_{r_N}^{[N]}.
\end{equation}

Crucially, this scheme does not require the evaluation of any mixed derivatives of $\qoi$ but only of its gradient, due to the fact that the operators $\U_{\dd_n}$ do not include tensor-products of osculatory interpolants;
thus, the derivative-induced curse of dimensionality of the tensor-product osculatory interpolants is broken.
Note that since in the setting of this section we only have two univariate levels,
convergence cannot be reached by enlarging the multi-index set,
that is fixed as $\indexsetsumw{1}$; rather, it can be achieved by increasing
the number of interpolation points,
\ie $r_n \rightarrow \infty$ along every direction $\uparVar_n$.

\begin{figure}[!t]
  \centering
  \includegraphics{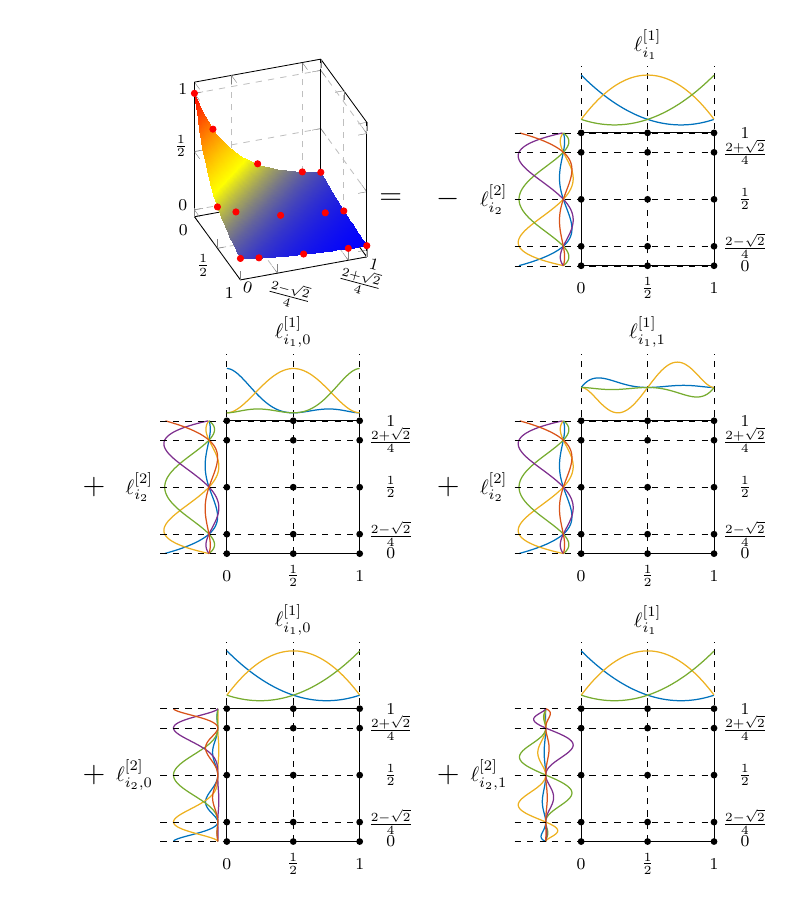}
  \caption{Osculatory interpolation by Smolyak algorithm on Cartesian grid with
    $3\times 5$ Clenshaw-Curtis points ($\bm r = [2,\,3]$ with $m$ doubling).
    The top-left panel shows $\SO_{\bm r}\qoi$ and the values of $\qoi$
    at the Cartesian grid $\interpPoints_{\bm r}$;
    the other 5 panels show $\interpPoints_{\bm r}$ and the combinations
    of nodal basis functions $\lagrange_{j}^{[n]}$, $\lagrange_{j,0}^{[n]}$,
    and $\lagrange_{j,1}^{[n]}$ used to build $\SO_{\bm r}\qoi$.
  }\label{fig:gh_1grid}
\end{figure}

Analyzing the polynomial spaces corresponding to each component operator
\[
\begin{aligned}
 \polyop{\U_{\bm 1}}&=\polyop{\L_{\bm r}}=\spn\Big\{\uparVar_1^{j_1}\ldots\uparVar_N^{j_N}: j_n\le m(r_n)-1, n=1\dots,N\Big\},
\\\polyop{\U_{\dd_n}}&=\spn\Big\{\uparVar_1^{j_1}\ldots\uparVar_N^{j_N}: j_n\le 2m(r_n)-1,\, j_s\le m(r_s)-1,\,  s=1,\dots,N,\, s\ne n \Big\},
\end{aligned}
\]
and using \cref{eq:polysparse} we conclude that
$\SO_{\bm r}$ delivers an approximation in the space
\begin{equation}\label{eq:single-grid-poly-space}
\polyop{\SO_{\bm r}}=\spn\Big\{\uparVar_1^{j_1}\ldots\uparVar_N^{j_N} : 0\le j_n\le 2m(r_n) -1, \text{ and }\#\{n: j_n \geq m(r_n)\}\le 1  \Big\}.
\end{equation}
To give further intuition on the construction of $\SO_{\bm r}$,
in \cref{fig:gh_1grid} we fix $N=2$, in which case $\SO_{\bm r}$ simplifies to
$\SO_{\bm r}= - \L_{r_1}^{[1]} \otimes \L_{r_2}^{[2]} + \L_{r_1}^{[1]} \otimes \O_{r_{2}}^{[2]} + \O_{r_1}^{[1]} \otimes \L_{r_{2}}^{[2]}$.  
The top-left panel shows the approximation $\SO_{\bm r}\qoi$ of 
$\qoi(\parVar) = (1 + \uparVar_1 + \uparVar_2)^{-3}$.
The remaining five panels show the combinations of nodal basis functions employed
to build $\SO_{\bm r}\qoi$, as obtained by expanding $\O_{r_{n}}^{[n]}$
in terms of $\lagrange_{j,0}^{[n]}$ and $\lagrange_{j,1}^{[n]}$ as in \cref{eq:uni-interp-osc}.

\subsection{Challenges in nesting the two Smolyak algorithms}\label{sec:challenges}
To solve both the derivative-induced and sampling curses of
dimensionality, \cite{deBaar2015} proposes to apply the Smolyak algorithm twice:
first as in \cref{sec:osculatory-smolyak}, constructing a
osculatory interpolant on a Cartesian grid using cubic splines
(\ie taking $\U_{1}^{[n]}$ as a $C^0$ piecewise linear interpolant and
$\U_{2}^{[n]}$ as the $C^1$ osculatory cubic spline interpolant at the same points);
second as in \cref{sec:sparse-grids}, combining several osculatory
spline interpolants on each component grid.
Using the same approach in our polynomial setting results in the following Smolyak operator 
\begin{equation}\label{eq:debaar-combitec}
  \SB_{w} \coloneqq \sum_{\bm i \in \indexsetsum}c_{\bm i}\SO_{\bm i},
\end{equation}
which delivers an approximation in the space
\begin{equation*}\label{eq:poly-broken}
  \polyop{\SB_w} = \sum_{\bm i \in \indexsetsum} \polyop{\SO_{\bm i}}.
\end{equation*}
Although tempting, this approach is however flawed, 
since the component operators $\SO_{\bm i}$ do not have a tensor-product
structure, therefore they do not admit a telescopic decomposition in terms of
detail operators as in \cref{eq:telescopic}, which is the cornerstone upon which
any combination technique is ultimately founded.

To investigate the implications of this fact we take a closer look to
$\SB_{w}$ as in \cref{eq:debaar-combitec}
and $\SL_{w}$ as in \cref{eq:lagr-sp-grid} built on the same sparse grid
$\interpPoints_{\indexsetsum}$
choosing $w=2$, Clenshaw--Curtis points and $m$ doubling.
For both approaches
the nonzero coefficients of the combination techniques  are $c_{(1,2)}$,
$c_{(1,3)}$, $c_{(2,1)}$, $c_{(2,2)}$, and $c_{(3,1)}$.

\begin{figure}[!t]
  \centering
  \includegraphics{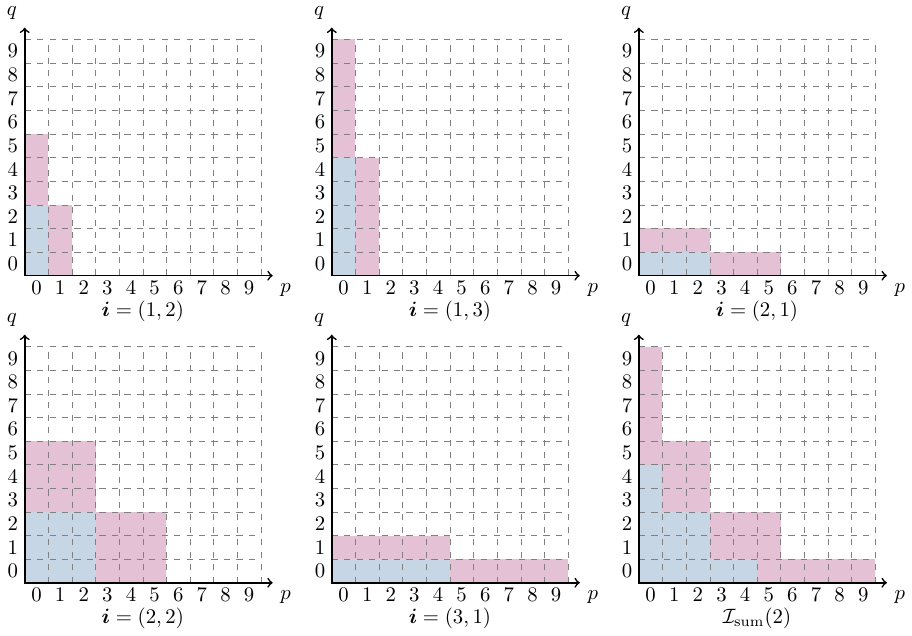}
  \caption{Polynomial spaces associated to each component grid 
    and to the resulting sparse grid: blue region for
    $\polyop{\L_{\bm i}}$ and $\polyop{\SL_{2}}$,
    union of the blue and magenta regions for $\polyop{\SO_{\bm i}}$
    and  $\polyop{\SB_{2}}$.}\label{fig:failure3}
\end{figure}

First, we observe that the dimension of $\polyop{\SB_2}$ is smaller than the
number of data (values and derivatives of $\qoi$) used by $\SB_{2}$,
implying that $\SB_{2}$ is not an interpolant even when employing nested points,
whereas $\SL_{2}$ is.
This argument is more easily followed with the help of \cref{fig:sg-plot} and
\cref{fig:failure3}. Specifically, \cref{fig:sg-plot} shows the five component
grids $\mathcal{T}_{\bm i}$ and the final sparse grid
$\interpPoints_{\indexsetsumw{2}}$. 
The sparse grid $\interpPoints_{\indexsetsumw{2}}$ has $13$ points, meaning that
$\SL_{2}$ employs 13 data (13 evaluations of $\qoi$) while $\SB_{2}$ employs $13
\times 3 = 39$ data (13 evaluations of $\qoi$ and $26=13 \times 2$ partial
derivatives). For each component grid $\mathcal{T}_{\bm i}$, the panel at the
corresponding location in \cref{fig:failure3} represents graphically the
ranges of the respective component operator of $\SL_{2}$ and $\SB_{2}$, \ie
polynomial spaces $\polyop{\L_{\bm i}}$ and $\polyop{\SO_{\bm i}}$
see equations \cref{eq:sparse-grid-poly-space} and
\cref{eq:single-grid-poly-space}. 
More specifically, each of these
panels shows a grid in which the square in position $(p,q)$ corresponds to the
monomial $\uparVar_1^p\uparVar_2^q$: the blue region thus indicates the polynomial
space $\polyop{\L_{\bm i}}$, while the union of the blue and magenta regions indicates
the polynomial space $\polyop{\SO_{\bm i}}$, where the magenta part indicates the
monomials brought in by having gradient information at disposal. The
bottom-right panel of the figure shows the final polynomial spaces
$\polyop{\SL_{2}}$ and $\polyop{\SB_2}$, that are obtained as the sum of the
component spaces. The dimension of the space and number of data for each
component grid and for the final approximations are reported in
\cref{tab:failure2}. The crucial observation is now that the dimension
of $\polyop{\SB_2}$ is 35 even though $\SB_{2}$ uses $39$ data.

\begin{table}[!t]\centering
  \footnotesize
  \begin{tabular}{ccccc}
    \toprule 
    & \multicolumn{2}{c}{$\SL_{2}$} & \multicolumn{2}{c}{$\SB_{2}$} \\
    \cmidrule(l){2-3} \cmidrule(l){4-5}
    grid	& space dim. & nb. data & space dim. & nb. data \\
    \midrule
    $(1,2)$ 	& 3 & 3 & 9 & 9 \\ 
    $(1,3)$ 	& 5 & 5 & 15 & 15 \\
    $(2,1)$ 	& 3 & 3 & 9 & 9 \\
    $(2,2)$ 	& 9 & 9 & 27 & 27 \\
    $(3,1)$ 	& 5 & 5 & 15 & 15 \\
    \midrule
    final count & 13 & 13 & 35 & 39 \\
    \bottomrule
  \end{tabular}\caption{Dimension of the polynomial space and number of data for each component grid and for the final approximation
    for $\SL_{2}$ and $\SB_{2}$.}\label{tab:failure2}
\end{table}

\begin{figure}[!t]
  \centering
  \includegraphics{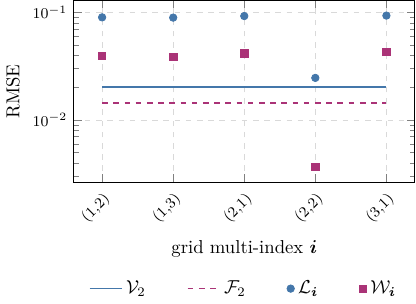}
  \caption{ RMSE for $\qoi(\parVar) = (1 + \uparVar_1 + \uparVar_2)^{-3}$
    when using the operators
    $\SL_{2}$, $\SB_{2}$ (horizontal lines) and their components $\L_{\bm i}$
    and $\SO_{\vec i}$ (markers), respectively.
  }\label{fig:failure1}
\end{figure}

Second, we compare the root mean square error based on $M_{mc}=1000$ random points
for the approximations $\SL_{2}\qoi$, $\SB_{2}\qoi$ of the function
$\qoi(\parVar) = (1 + \uparVar_1 + \uparVar_2)^{-3}$
and their component operators, \ie
\begin{equation}\label{eq:qoi-rmse}
  \mathrm{RMSE} \coloneqq \sqrt{\frac{1}{M_{mc}} \sum_{j=1}^{M_{mc}}\left[\qoi(\parVar_j)-\aqoi(\parVar_j)\right]^2},
\end{equation}
where $\aqoi$ is a placeholder for $\SL_{2}\qoi$, $\SB_{2}\qoi$, and $\L_{\bm
  i}\qoi$, $\SO_{\bm i}\qoi$ for $\bm i \in \indexsetsumw{2}$. The RMSE values
are reported in \cref{fig:failure1}: the RMSE of the final approximations
$\SL_{2}\qoi$ and $\SB_{2}\qoi$ are indicated by horizontal lines whereas the
RMSE for $\L_{\bm i}\qoi$ and $\SO_{\bm i}\qoi$ are marked by circles and
squares, respectively, and are arranged in columns corresponding to their
multi-index $\bm i$. As expected, $\SO_{\bm i}\qoi$ is more accurate than
$\L_{\bm i}\qoi$ for all $\bm i$. However, while the RMSE of $\SL_{2}$ is
\emph{smaller} than that of its components (the blue solid line is \emph{below}
the blue circles), the RMSE of $\SB_{2}$ is \emph{larger} than that of its
component $\SO_{(2,2)}\qoi$ (the magenta dashed line is \emph{above} the magenta
square for $\bm  i=(2,2)$); in other words, a single component operator deliver
better approximation than $\SB_{2}$, making it evident that this approach shows
issues in practical applications.

Summarizing, the nested application of the Smolyak algorithms is theoretically
unjustified and does not provide neither interpolatory properties nor good
approximations. Note that in \cref{fig:failure1} $\SB_{2}$ seems to be
nonetheless competitive with $\SL_{2}$, and a similar situation occurred in
\cite{deBaar2015}, in which the piecewise cubic osculatory sparse grid was
observed to converge and even apparently improve upon the classical piecewise
linear sparse grid. 

\subsection{A possible remedy: least-squares on the sparse grid}\label{sec:GELS}

The discussion in the previous section naturally leads to replacing formulation
\cref{eq:debaar-combitec} with a least-squares approach. More specifically,
given a multi-index set $\indexsetsum$ and the evaluations of $\qoi$ and $\nabla
\qoi$ on the sparse grid $\interpPoints_{\indexsetsum}$, we compute the final
approximation $\G_{w}\qoi$ as the (weighted) least-squares approximation of
$\qoi$ on $\polyop{\SB_w}$:
\begin{equation}\label{eq:least-squares-approximant}
    \G_{w}\qoi \coloneqq \argmin_{v\in \polyop{\SB_w}} \sum_{\bm \interpPoint_{\bm j} \in \interpPoints_{\indexsetsum}} \left[ \left\|\qoi(\bm \interpPoint_{\bm j}) - v(\bm \interpPoint_{\bm j})\right\|^2 + \sum_{n=1}^N \omega_n\left\| \partial_{\uparVar_n} \qoi(\bm \interpPoint_{\bm j}) - \partial_{\uparVar_n} v(\bm \interpPoint_{\bm j})\right\|^2 \right],
\end{equation}
with suitable weights $\omega_n > 0$ that can be used to balance the contributions
of the function and gradient values when the their magnitudes are very different.
We call this method \emph{Gradient-Enhanced Least-squares on Sparse grids} (\gels).
In practice, we construct the set of Legendre orthonormal polynomials spanning $\polyop{\SB_w}$ (\ie a polynomial chaos expansion)
and use this basis to solve the minimization problem in \cref{eq:least-squares-approximant}; the derivatives of the Legendre polynomials are computed following \cite{barrio2002}. 
Convergence of $\G_{w}\qoi$ to $\qoi$ is then expected as the Smolyak level $w$ grows,
which increases the cardinality of the sparse grid $\interpPoints_{\indexsetsum}$ and correspondingly enlarges the polynomial space $\polyop{\SB_w}$.

In other words, we propose a hybrid, \emph{sparse-grid inspired} approach, in
which a Smolyak procedure defines a sparse-grid sampling scheme
and a suitable polynomial space for the approximation, but the coefficients of
such approximation are computed by a least-squares approach rather than by
interpolation. Note in particular that by using a minimization principle we are
enforcing that, at least in the discrete norm, the error of $\G_w\qoi$ is
smaller than that of $\SO_{\bm i}\qoi$ for any $\bm i \in \indexsetsum$, solving
the issue discussed in \cref{fig:failure1} (the RMSE of \gels would be $3.098\cdot 10^{-3}$ that is lower than the smallest error marked in the figure, \ie $3.65 \cdot 10^{-3}$ for $\SO_{\bm 2}Q$).
Note also that our approach differs from the procedures typically
employed in literature to build a polynomial chaos expansion
by a least-squares-based approach (either gradient-enhanced as the already-mentioned
\cite{guo.eal:2018,jakeman2015, peng2016} or not,
as in \eg \cite{sudret:adaptive.pce.with.reg,migliorati:rdpforfunc})
in at least two points:
first, the sparse-grid samples in $\interpPoints_{\indexsetsum}$ are not space-filling,
whereas this property is typically sought after;
second, the balance between the number of samples and the selection of the
polynomial basis is based neither on usual properties
such as the restricted isometry property of the design matrix
\cite{guo.eal:2018, jakeman2015, peng2016} or its condition number
\cite{migliorati:rdpforfunc},
nor on adaptive algorithms \cite{sudret:adaptive.pce.with.reg,jakeman2015}.
Theoretical considerations about these aspects of our approach
are certainly worth pursuing but are left to future works.

\section{Numerical results}\label{sec:results}

In the section we investigate the performance of \gels
comparing it against that of the standard sparse-grid method
\cref{eq:lagr-sp-grid} (\sg in the following)
in different scenarios of increasing complexity.
In \cref{sec:numeric-genz} we consider a curated
selection of benchmark functions  $\qoi$  \emph{known in closed-form} having
different regularity and different domain dimension $N$. In this case we also
compare our strategy against a version of \cref{eq:least-squares-approximant}
where $\interpPoints_{\indexsetsum}$ is replaced by Monte Carlo samples, to
verify the benefits of using sparse-grid sampling in such context. In
\cref{sec:numeric-noise} we still consider functions in closed-form, that we
however artificially pollute with random noise to investigate the
\emph{robustness} of \gels when working with \emph{data that are known with
different levels of accuracy}. This situation is the standard scenario when
working with parametric PDE, since we can expect that typically \QoI and their
gradients are known with different precision. In \cref{numeric-plap}, we finally
move to the actual use-case for which our methodology is intended, \ie computing
a surrogate model for a \emph{parametric PDE}. In detail, we consider two
versions of a non-linear Darcy flow problem ($p$-Laplace equation): we show that
only in the first case the gradient enhancement improves the accuracy with
respect to the standard sparse grids, and close the section providing a
heuristic explanation for this behavior.

In all these scenarios,
having checked that $\qoi$ and $\partial_{\uparVar_n}\qoi$ have the same 
order of magnitudes in all experiments,
we set the weights $\omega_n$
in \cref{eq:least-squares-approximant} to $1$, and  
monitor the decay of the approximation error as a function of the computational cost:
\begin{itemize}
\item the \emph{error} is measured as the RMSE \cref{eq:qoi-rmse}
  between the quantity of interest
  and the different approximations of $\qoi$ considered in each scenario.
\item the \emph{cost} is defined solely in terms of function and gradient evaluations, \ie neglecting \eg
  the costs of computing the collocation points or solving the least-squares system,
  as well as the factorization issues discussed in Remark~\ref{rem:assembly-costs}.
  We further assume that the cost of evaluating $\qoi(\parVar)$ is constant for every $\parVar \in \parSpace$
  (namely, that the number of Newton iterations does not depend on $\parVar$), and take such cost as a cost unit.
  Finally, we denote by $\lambda$ the relative cost of computing a single partial derivative $\partial_{\uparVar_n}\qoi$:
  owing to \cref{eq:par-derivative-of-qoi} and following discussion, we have that $0\leq \lambda \leq 1$, and 
  we can quantify the cost of the methods as follows
  \begin{align}
    \text{cost per point with \sg} & : 1, \nonumber\\
    \text{cost per point with \gels} & : 1 + N \lambda, \ \text{with}\ \lambda\in[0,1].\label{eq:costHSG}
  \end{align} 
\item the sequence of approximations for which cost and error are computed is obtained by increasing the Smolyak level $w=1,2,\ldots$ in the construction of \gels and \sg.
\end{itemize}
The results 
were obtained using the Sparse Grids Matlab Kit \cite{piazzola2024}.

\subsection{Benchmark functions}\label{sec:numeric-genz}

\begin{figure}[!t]
	\centering
	\includegraphics{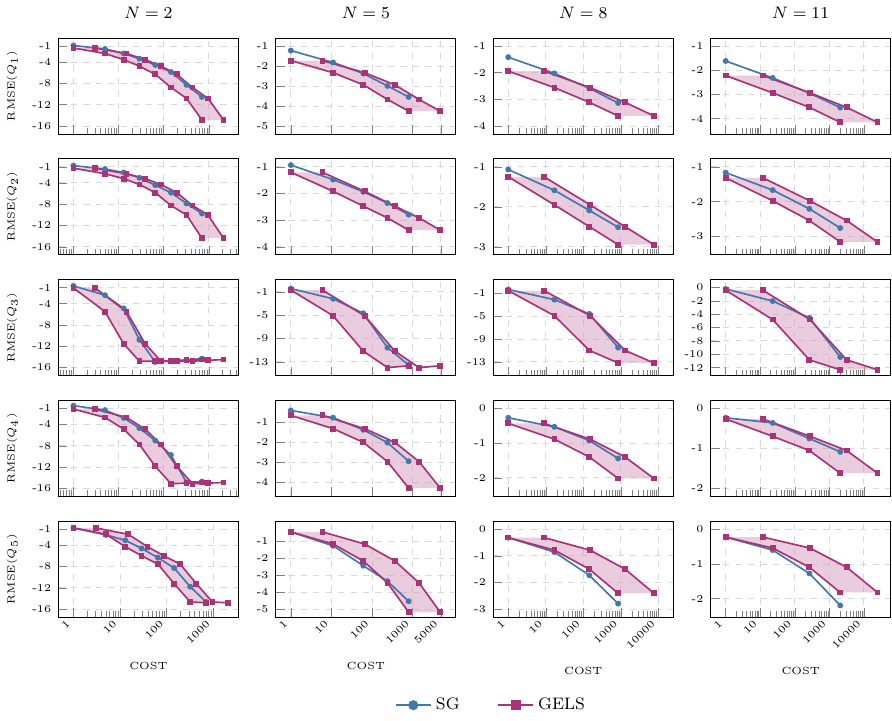}
	\caption{
        $\log_{10}(\text{RMSE})$ as function of cost for both
        \sg (lines) and \gels (strips corresponding to $\lambda\in[0,1]$)
        for the \QoI in \cref{eq:gn}.
        The scale of the horizontal axis is common in each column;
        the scale of the vertical axis is different for every plot.
      }
      \label{fig:exm-all-genz}
\end{figure}

In this section we consider a selection of benchmark functions with different behaviors
on $\parSpace = [0,1]^N$ for $N=2, 5, 8, 11$, see \eg \cite{barthelmann.novak.ritter:high},
namely
\begin{equation}\label{eq:gn}
\begin{aligned}
	\qoi_1(\parVar) & = \frac{1}{1+\sum_{n=1}^N\uparVar_n}, &
	\qoi_2(\parVar) & = \frac{1}{1+\sum_{n=1}^N\uparVar^2_n}, &
	\qoi_3(\parVar) & = \sum_{n=1}^N\cos(\uparVar_n), \\
	\qoi_4(\parVar) & = \cos\Big(2\pi + \sum_{n=1}^N \uparVar_n\Big), & 
	\qoi_5(\parVar) & = \exp\Big(-\sum_{n=1}^N \big(\uparVar_n - \frac{1}{2}\big)^2\Big),&
\end{aligned}
\end{equation}
and report in  \cref{fig:exm-all-genz} the convergence of the RMSE of the \sg and \gels approximations.
To account for the variability of cost of \gels due to $\lambda$ in \cref{eq:costHSG},
we actually plot the two convergence curves corresponding
to the limit cases $\lambda=0$ (each partial derivative is obtained at zero cost) and $\lambda=1$ (each partial derivative is as expensive as a function evaluation)
and shade the entire strip between them.
This is different from what done in the papers about gradient-enhanced approximation
mentioned in the introduction,
that typically assume a specific value for $\lambda$ when reporting convergence as a function of cost.

\begin{figure}[!t]
	\centering
	\includegraphics{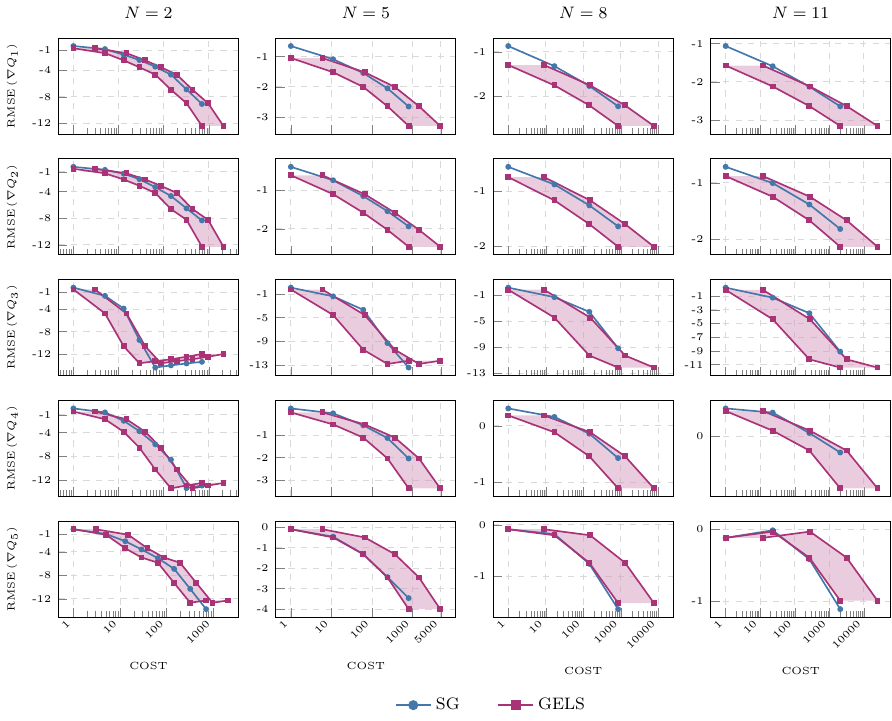}
    \caption{
        $\log_{10}(\text{RMSE})$ as function of cost for both
        \sg (lines) and \gels (strips corresponding to $\lambda\in[0,1]$)
        for the partial derivatives of \QoI in \cref{eq:gn}.
        The scale of the horizontal axis is common in each column;
        the scale of the vertical axis is different for every plot.
	}
	\label{fig:exm-all-genz-grad}
\end{figure}

We can observe that for $\qoi_1$, $\qoi_3$, $\qoi_4$
\gels outperforms \sg for a large range of values of $\lambda$.
For the benchmark $\qoi_2$ 
the range of $\lambda$ for which \gels outperforms \sg gets instead smaller and smaller as $N$ increases,
and the situation is even worse for the benchmark $\qoi_5$,
in which \gels outperforms \sg only for $N=2$ and if the derivatives can be computed at a limited cost.
Similar conclusions on the \sg and \gels accuracy can be drawn also by looking at the trend of the RMSE of the gradient of $\qoi$ (more precisely, at the square root of the sum of the mean square errors of each partial derivative), see \cref{fig:exm-all-genz-grad}.
This suggests that the usefulness of enhancing a sparse-grid-like approximation with gradient information is problem-dependent
(both in terms of regularity of $\qoi$ with respect to $\parVar$ and of cost of evaluating derivatives of $\qoi$),
and a fine-tuned theoretical analysis should be carried out to understand the approximation power of \gels.

As a second experiment, we compare the convergence of \gels against
gradient-enhanced least squares approximation based on Monte-Carlo (MC) samples
in the same polynomial spaces $\polyop{\G_w}$ and with same number of sample points
as $\interpPoints_{\indexsetsum}$.
The aim of this test is to verify whether sampling $\qoi$ and its gradients
over a sparse grid pays off, at least pre-asymptotically,
with respect to a standard choice for least-squares such as MC.
Note that besides a sampling scheme our approach gives a practical recipe to design the polynomial space where the least-squares approximation is sought
and the collocation points are deterministic, whereas using a MC-based approach leaves the user with the question of
how to design a suitable polynomial space and with a result that is sample-dependent: in this sense,
even a result in which the MC-based approach is comparable/slightly better would still be an acceptable outcome.

\begin{figure}[!t]
	\centering
	\includegraphics{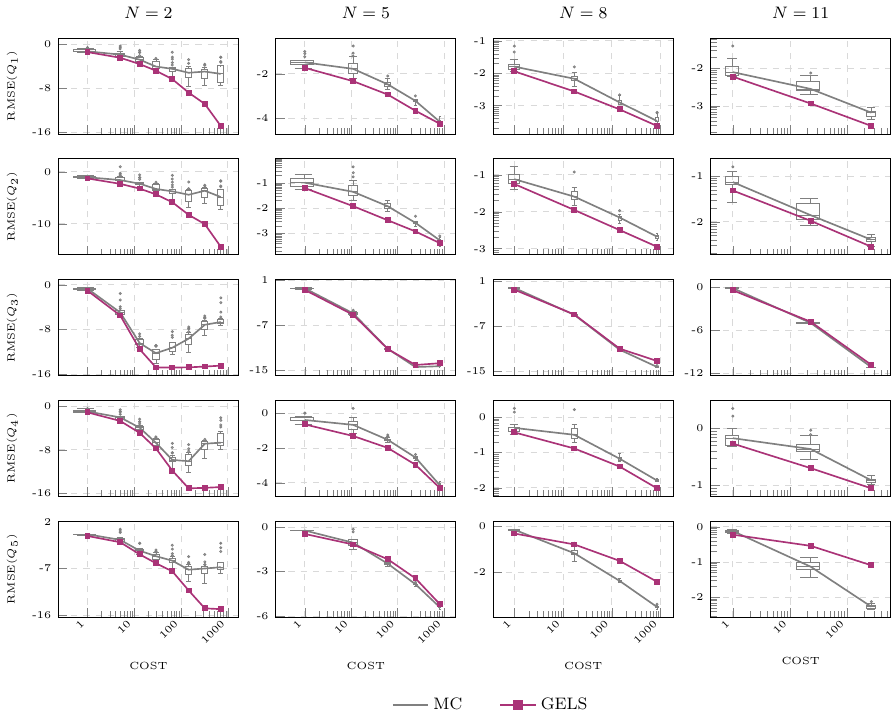}
	\caption{
        $\log_{10}(\text{RMSE})$ as function of cost for \gels and MC-based least-squares for the \QoI in \cref{eq:gn}.
        The scale of the horizontal axis is common in each column;
        the scale of the vertical axis is different for every plot.
        }
	\label{fig:exm-mc-rmse}
\end{figure}

The results of this comparison are shown in \cref{fig:exm-mc-rmse}. As for the previous \cref{fig:exm-all-genz,fig:exm-all-genz-grad},
different rows correspond to the different test functions in
\cref{eq:gn}
and different columns to different parameter space dimensions $N=2,5,8,11$. 
Since both methods employ values and partial derivatives of $\qoi$,
they have the same computational cost per point, that we take as unitary cost.
For MC, we repeat the analysis 30 times, and report in the plot
statistics of the errors across the repetitions, namely the
convergence of the median error and the boxplot with outliers
for each number of collocation points.
We can observe that for $N=2$ the error of MC stagnates or even diverges, which is the typical behavior
whenever the number of collocation points is not enough to guarantee that the least-squares system has a good condition number,
see \eg \cite{migliorati:rdpforfunc}.
Conversely, \gels seems to be robust with respect to this issue and delivers a better convergence rate in the range of cost for which MC does not diverge.
The onset of this phenomenon is not yet visible within the range of computational costs considered for larger values of $N$;
however, \gels nonetheless shows a better performance than MC-based-least-squares for
$\qoi_1, \qoi_2, \qoi_4$,
a comparable performance for $\qoi_3$, and worse performance only for $\qoi_5$,
for which even the previous test was showing a poor convergence of \gels.
Furthermore, also in this test similar conclusions can be drawn by looking at the convergences of the RMSE of the approximation of the partial derivatives of $\qoi$ (not shown).
Thus, we can conclude that our choice of using a sparse grid as set of collocation points is sensible and works well for the designed polynomial space.

\subsection{Data with different accuracy}\label{sec:numeric-noise}

\begin{figure}[!t]
  \centering
  \includegraphics{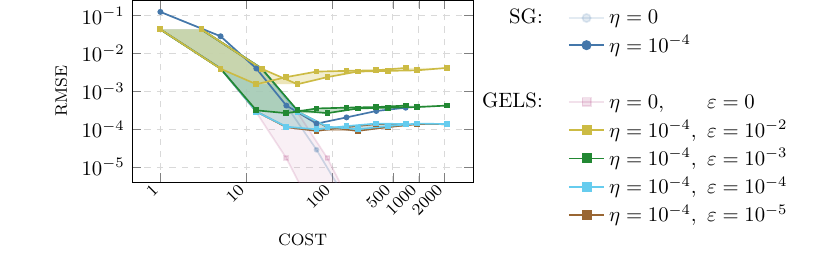}
  \caption{
  	RMSE as a function of the computational cost for both
  \sg (lines) and \gels (strips corresponding to $\lambda\in[0,1]$)
  depending on the accuracy of the values, perturbed by a normal noise with standard deviation $\eta$, and the gradient information, perturbed by a normal noise with standard deviation $\varepsilon$.
  }
\label{fig:noise}
\end{figure}

The examples presented in the previous section focus on the approximation of \QoI whose values and gradients can be evaluated exactly; 
nevertheless, in practical applications, such data are typically known with different precisions and/or noisy, due to measurement or numerical errors.  
Therefore, in this section, we examine the performance of \gels
with respect to the relative accuracy with which values and gradients of $\qoi$ are known:
as can be expected, this relative accuracy can slow down or even harm the convergence of \gels with respect to its gradient-less counterpart \sg
(some discussion in this sense is also reported in \cite{bhaduri2020,deBaar2015,peng2016}).

In this test we thus consider $\qoi_1$ defined in \cref{eq:gn} on $\parSpace = [0,1]^2$,
with $N=2$ and pollute its value and gradient evaluations with random noise
of different amplitude. 
More precisely, we introduce the perturbations $Z_\eta(\parVar)$,
$Z_\varepsilon(\parVar)$ that for each $\parVar$
are realizations of Gaussian random variables with zero mean and standard deviations
$\eta$ and $\varepsilon$, respectively, and denote by
$\SL_{w,\eta}\qoi$ the \sg operator where the evaluations $\qoi(\parVar)$
are replaced by $\qoi(\parVar) + Z_\eta(\parVar)$, and by
$\G_{w,\eta,\varepsilon}\qoi$ the \gels operator where both the evaluations
$\qoi(\parVar)$ and $\partial_{\uparVar_n}\qoi(\parVar)$ are replaced by
$\qoi(\parVar) + Z_\eta(\parVar)$ and $\partial_{\uparVar_n}\qoi(\parVar) +
Z_\varepsilon(\parVar)$.

We then consider several combinations of $\eta$ and $\varepsilon$ and show the
resulting error convergence of the perturbed \sg and \gels operators in
\cref{fig:noise}; the RMSE are 
computed against the \emph{unperturbed} values of $\qoi$. 
In detail, we plot the convergence curve/strip for:
\begin{itemize}
\item the unperturbed approximations $\SL_{w}\qoi=\SL_{w,0}\qoi$ and
  $\G_{w}\qoi=\G_{w,0,0}\qoi$, that serve as the reference/benchmark
  convergences, \ie they mark the best achievable accuracies for the two
  methods: as already reported in the top left panel with
  \cref{fig:exm-all-genz} this test is such that exact knowledge of
  $\partial_{\uparVar_n} \qoi$ is beneficial for approximating $\qoi$;
\item the perturbed \sg approximation $\SL_{w,10^{-4}}\qoi$; as can be
  expected, this approximation converges identically to its unperturbed
  counterpart as long as the approximation error is larger than the noise level,
  after which the convergence stagnates;
\item the perturbed \gels approximations $\G_{w,10^{-4},\varepsilon}\qoi$ for
  $\varepsilon= 10^{-5}$, $10^{-4}$,  $10^{-3}$,  $10^{-2}$; we can
  observe that, as expected, the convergences of $\G_{w,\eta,\varepsilon}\qoi$
  follow the unperturbed counterpart until they 
  stagnate at an error of about $\max\{\eta,\varepsilon\}$.
\end{itemize}
In conclusion, this test shows that the incorporating evaluations of $\partial_{\uparVar_n} \qoi$
improves the approximation not only if these evaluations can be obtained for a sufficiently small cost,
but also if their accuracy is comparable to that of the evaluations of $\qoi$, which is not necessarily true \eg when $\qoi$ and $\partial_{\uparVar_n} \qoi$ are being obtained
from numerically solving the respective PDEs over the same mesh.
The next section further elaborates on this topic.

\subsection{An example with a parametric PDE}\label{numeric-plap}

\begin{figure}[!t]
  \centering
  \includegraphics{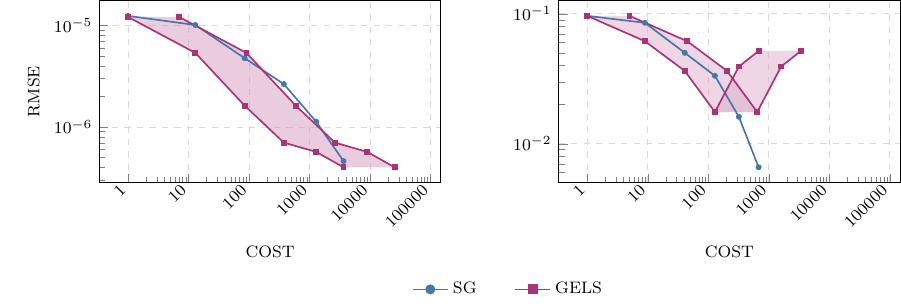}
  \caption{
  Trend of the RMSE against the computational cost for \sg (lines) and \gels (strips corresponding to $\lambda\in[0,1]$).
  Left: test1, right: test2.
}
\label{fig:pde-test}
\end{figure}

In this section, we compare the performance of \gels and \sg on two test cases in which evaluating the quantity of interest $\qoi$ requires solving a non-linear parametric PDE,
following the setup of \cref{sec:parametricPDEs}.
In the first case \gels will be better than \sg, whereas in the second case the results are the opposite, see \cref{fig:pde-test}:
in the rest of the section we explain this difference in view of the experience gathered in the previous sections.

In detail, we consider a Darcy-type flow with heterogeneous permeability, modeled as a $p$-Laplace equation defined on a planar domain $\Omega=[-1,1]^2$,
\begin{equation}\label{eq:plap}
  \begin{cases}
    - \mathrm{div}\left(K(\parVar) \,\vert \nabla \sol\vert^{p-2} \nabla \sol\right) = f, &\text{in}\ \Omega\\
    \sol = g, &\text{on}\ \partial\Omega
  \end{cases}
\end{equation}
where the data $f$, $g$, the permeability field $K(\parVar)$, and
the parameter space $\parSpace$, are specified in each test case.
In what follows we consider $p=1.8$, that is in the range $3/2\le p\le 2$
usually considered for nonlinear filtration in porous media  and ice-sheet modeling
\cite{benedikt2018,Jakeman2025}. 
The quantity of interest for this problem is the total flux through the right boundary of the domain,
$\Gamma_R = \{1\}\times[-1,1]$, \ie
\begin{equation*}
  \qoi(\parVar) = \int_{\Gamma_R} ( - K(\parVar) | \nabla \sol |^{p-2} \nabla \sol ) \cdot \bm n\textrm{ d}\Gamma_R.
\end{equation*}
Equation \cref{eq:plap} is solved numerically using a Galerkin method based on tensor-product $C^2$ cubic spline basis functions on uniform meshes.
Note that $\qoi(\parVar)$ depends on $\parVar$ not only indirectly through $\sol$ but also
directly through $K(\parVar)$, therefore a suitable modification of \cref{eq:par-derivative-of-qoi} is needed to compute $\partial_{\uparVar_n}\qoi$. 

The RMSE of \gels and \sg in these tests are computed not against exact values
of $\qoi$ (that are clearly unavailable), but rather against
approximated values obtained solving \cref{eq:plap}
with a mesh with  $450\times450$ elements.

\subsubsection*{Test 1 -\gels better than \sg}

\begin{figure}[!t]
  \centering
  \includegraphics{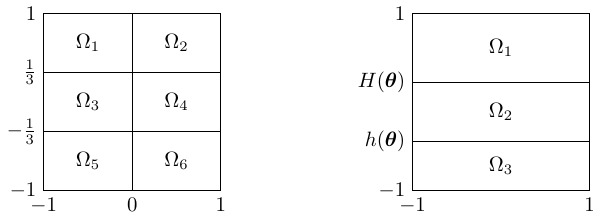}
  \caption{Left: partition of $\Omega$ for Test 1. Right: partition of $\Omega$ for Test 2.}\label{fig:plap-domains}
\end{figure}

In this test, we set $f=0$ and choose a Dirichlet boundary condition
that decreases linearly from $g=2$ on the left-hand boundary to $g=1$ on the right-hand boundary, \ie $g(\bm x) =(3 - x_1)/2$.
The permeability field is a piecewise-constant function over a $3\times2$ chessboard-like partition of $\Omega$,
see \cref{fig:plap-domains}-left, 
with each value being controlled by a different parameter $\uparVar_n$:
\begin{equation*}
  K(\bm x;\parVar) = \sum_{n=1}^6 10^{\uparVar_n} \chi_n(\bm x) , \quad \bm x \in \Omega,
\end{equation*}
where $\chi_n$ are the indicator functions of the six sub-regions $\Omega_n$ numbered as in \cref{fig:plap-domains}-left, and $\parVar \in [-9,-4]^N$.
To cope with the discontinuities of $K$, the basis functions of the Galerkin method
are modified to be $C^0$ along the discontinuity lines.

\begin{figure}[!t]
  \centering
  \includegraphics{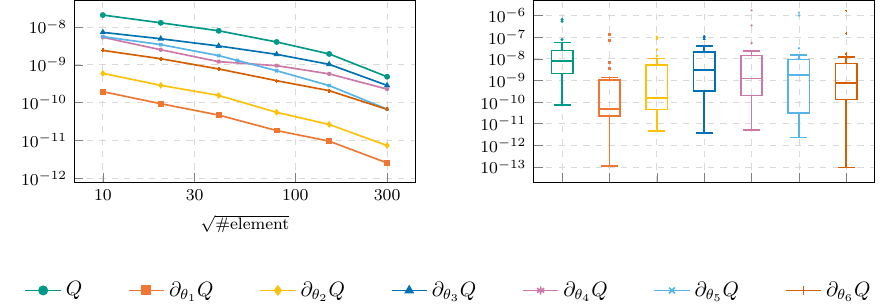}
  \caption{
  Test 1.
  Left:
  medians (over 20 realizations of $\parVar$) of the absolute error of the numerical approximation of
  $\qoi$ and $\partial_{\uparVar_n}\qoi$ on increasingly refined meshes.
  Right: Boxplot of the same numerical errors on the $40 \times 40$ mesh.
}
  \label{fig:tikz_3x2chess_f0_boxplot_convergence}
\end{figure}

Since Dirichlet data are prescribed on the entire boundary and the source term is set to zero,
the maximum principle ensures that all extrema of $\sol$ occur on $\partial\Omega$.
Consequently, the solution and its parametric derivatives remain well-behaved 
in the interior of the domain, with no localized regions of steep gradients or singular behavior.
We can therefore expect that the numerical approximation of $\qoi$ and of $\partial_{\uparVar_n}\qoi$ will converge with the same rate with respect to the number of
  elements of the mesh when approximating them by a Galerkin scheme.
To confirm this intuition, we compute the absolute error of the numerical
approximation of $\qoi$ and $\partial_{\uparVar_n}\qoi$ for a set of meshes of $\Omega$
with increasing number of elements and for 20 random realizations of
$\parVar \in \parSpace$: \cref{fig:tikz_3x2chess_f0_boxplot_convergence}-left
shows the convergence of the medians of such errors over each mesh, and indeed
they all appear to be converging with the same rate as the meshes gets finer.
More importantly, the medians of the errors for $\qoi$ are larger than
those for $\partial_{\uparVar_n}\qoi$: this suggests that  
for this test the derivative information is at least as reliable as the \QoI itself,
and hence we can expect the convergence of \gels to be better than that of \sg
(if the cost of computing $\partial_{\uparVar_n}\qoi$ is small enough).
Further confirmation of the comparison among the errors on $\qoi$
and $\partial_{\uparVar_n}\qoi$ is obtained from
\cref{fig:tikz_3x2chess_f0_boxplot_convergence}-right, that shows the
boxplots of the errors for the 20 random realizations of $\parVar$ on the
$40\times 40$ mesh.

Upon fixing the mesh to $40\times 40$ elements, we perform the usual convergence
study for \gels and \sg.
The convergence of the two methods is shown in \cref{fig:pde-test}-left: 
the convergence of \gels is faster than that of \sg until the error reaches $10^{-6}$,
at which level the Galerkin approximation of $\qoi$ and that of $\partial_{\uparVar_n}\qoi$
may become comparable for some $\parVar$, see again
\cref{fig:tikz_3x2chess_f0_boxplot_convergence}-right.

\subsubsection*{Test 2 - \gels worse than \sg}

  In this second test we consider a continuous, piecewise linear permeability field $K$ stratified into three horizontal layers. 
  Contrary to the previous test however, the position of the three layers is not fixed, but is also parametrized by $\parVar$,
  which significantly increases the complexity of the function $\qoi(\parVar)$ and of its derivatives.

Specifically for $\parVar\in\parSpace:=[0,1]^4$, the domain $\Omega$ is divided in three subdomain as
in \cref{fig:plap-domains}-right, where the interfaces are at
$h(\parVar)=0.8\,\uparVar_3 - 0.9,$ and $H(\parVar)=0.8\,\uparVar_4 + 0.1$,
and the permeability $K$ is set to
\[
  K(\bm x;\parVar) = \sum_{m=1}^3 K_m(\parVar) \chi_m(\bm x;\parVar),
\]
where $\chi_m$ is the indicator function of $\Omega_m$ and 
$K_1(\parVar)= 10^{5\,\uparVar_2-9}$,
$K_3(\parVar)= 10^{5\,\uparVar_1-9}$,
$K_2(\parVar)=[K_3(\parVar)\left( H(\parVar)-x_2\right)+K_1(\parVar)\big( x_2-h(\parVar)\big)]/[H(\parVar)-h(\parVar)]
$.
The subdomain interfaces range in the intervals $[-0.9, -0.1]$ and $[0.1, 0.9]$, respectively, so that $\Omega_2$ is never empty. The permeabilities $K_1$, $K_2$ and $K_3$
range within $[10^{-9}, 10^{-4}]$ as in the previous test. 
We complete the formulation of \cref{eq:plap} by setting homogeneous Dirichlet boundary conditions and $f=1$.

\begin{figure}[!t]
	\centering
	\includegraphics{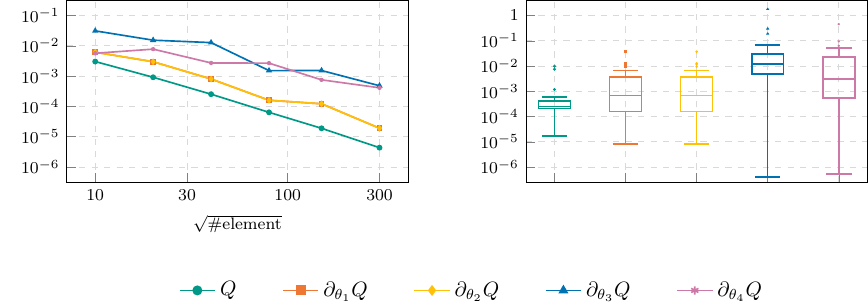}
  \caption{
	Test 2.
	Left:
	medians (over 20 realizations of $\parVar$) of the absolute error of the numerical approximation of
	$\qoi$ and $\partial_{\uparVar_n}\qoi$ on increasingly refined meshes.
	Right: Boxplot of the same numerical errors on the $40 \times 40$ mesh.
}
	\label{fig:tikz_3moving_layers_boxplot_convergence}
\end{figure}

This setup has an important structural consequence: the maximum principle no longer applies,
and the solution develops extrema in the interior of $\Omega$.
As a result, $|\nabla \sol|$ vanishes inside $\Omega$,
producing regions where the non-linear diffusion
term $K(\parVar)|\nabla \sol|^{p-2}\nabla \sol$ becomes degenerate.
In turn, this leads to a loss of spatial regularity of the solution
that further exhacerbates in its parametric derivatives,
impacting the accuracy of the respective spatial numerical approximations.
As done in the previous test, we compute the absolute error of the numerical
approximation of $\qoi$ and $\partial_{\uparVar_n}\qoi$ for a set of meshes of $\Omega$
with increasing number of elements and for 20 random realizations of
$\parVar \in \parSpace$. 
\cref{fig:tikz_3moving_layers_boxplot_convergence}-left shows the convergence of
the median errors of $\qoi$ and of $\partial_{\uparVar_n}\qoi$
as the spatial mesh is refined: while the error of $\qoi$ decreases steadily,
the errors of $\partial_{\uparVar_n}\qoi$ remain substantially larger,
typically 1 to 2 orders of magnitude larger than the error of $\qoi$ for the same mesh size.
This discrepancy is illustrated more clearly in 
\cref{fig:tikz_3moving_layers_boxplot_convergence}-right that shows the distribution across 20 random
realizations of $\parVar \in \parSpace$ on the mesh with  $40\times40$ elements, that 
will be used for studying the convergence of \sg and \gels.
\cref{fig:pde-test}-right then consistently shows that the \sg approximation
converges smoothly as computational cost increases,
whereas the error of \gels quickly stagnates.

\section{Conclusions}
\label{sec:conclusions}

This work discusses how incorporating gradient information in the construction of sparse-grids-like approximations
of multivariate functions is a challenging task.
From a theoretical point of view, a naive tensorization of osculatory
polynomials would run not only into the standard sampling curse of dimensionality but also into what we called
derivative-induced curse of dimensionality. Constraining ourselves to employing only gradient information
(without mixed derivatives) proves to break the usual combination technique, such that we resort to a hybrid approach combining sparse-grids and least-squares (\gels). The performance of \gels with respect to
the a standard sparse-grid approximation is then crucially depending on two factors: the cost of actually
obtaining the gradient information (which must be small enough to pay off), and the relative accuracy
of the evaluations of such gradients with respect to the function evaluation, which is not always granted
when both come \eg from the discretization of a PDE.

\section*{Acknowledgments}
      AB, FL, LT have been partially supported by the project PRIN PNRR ``Uncertainty Quantification of coupled models for water ﬂow and contaminant transport'' (P2022LXLYY), financed by the European Union -- NextGeneration EU.
      SI, LT have been partially supported by ICSC-Centro Nazionale di Ricerca in High Performance Computing, Big Data, and Quantum Computing, funded by European Union -- NextGeneration EU;
      AB, FL have been partially supported by the project COSMIC (PRIN 2022A79M75) of MUR funded by the European Union -- NextGeneration EU; 
      LT, AB, FL are members of the Gruppo Nazionale Calcolo Scientiﬁco-Istituto Nazionale di Alta Matematica (GNCS-INdAM).

\bibliographystyle{plain}
\bibliography{sample}

\end{document}

%% file: ex_shared.tex
\usepackage{amsmath}
\usepackage{mathrsfs}
\usepackage{amssymb}
\usepackage{amsfonts}
\usepackage{amsthm}

\usepackage{graphicx}
\usepackage{epstopdf}

\usepackage{multirow}
\usepackage{tabularx}
\usepackage{longtable}
\usepackage{booktabs} 
\usepackage{array}    

\usepackage{xspace}
\usepackage{bm}
\usepackage{cleveref}

\ifpdf
  \DeclareGraphicsExtensions{.eps,.pdf,.png,.jpg}
\else
  \DeclareGraphicsExtensions{.eps}
\fi

\usepackage{enumitem}
\setlist[enumerate]{leftmargin=.5in}
\setlist[itemize]{leftmargin=.5in}

\newtheorem{remark}{Remark}

\title{Sparse-grids-like surrogate models enhanced with gradient information}

\author{Andrea Bressan$^1$, Sofia Imperatore$^{1,2}$,\\
Francesca Locatelli$^1$, Lorenzo Tamellini$^1$
  }

\date{
\today\\[1ex]
\small
$^1${
    Istituto di Matematica Applicata e Tecnologie Informatiche ``E. Magenes'', Consiglio Nazionale delle Ricerche,
      Via Ferrata 5/A, Pavia, Italy 
  (\texttt{andrea.bressan@imati.cnr.it}, \texttt{francesca.locatelli@imati.cnr.it}, \texttt{lorenzo.tamellini@imati.cnr.it}).}
\\[1ex]
$^2${Institute for Complex Molecular Systems, Department of Biomedical Engineering, Eindhoven University of Technology, Eindhoven, The Netherlands
  (\texttt{s.imperatore@tue.nl}).}
}

\usepackage{amsopn}

\newcommand{\ie}{i.e.\ }

\newcommand{\eg}{e.g.,\xspace}

\newcommand{\coloneqq}{:=}

\input{notation}


%% file: notation.tex
\newcommand{\QoI}{QoI\xspace}

\renewcommand{\vec}[1]{\bm{#1}}

\newcommand{\N}{\mathbb{N}}
\newcommand{\Np}{\N_{+}}
\newcommand{\R}{\mathbb{R}}

\newcommand{\G}{\mathcal{G}} 
\newcommand{\U}{\mathcal{U}}  
\renewcommand{\L}{\mathcal{L}}
\renewcommand{\O}{\mathcal{O}}

\renewcommand{\S}{\mathcal{S}} 
\newcommand{\SL}{\mathcal{V}}  
\newcommand{\SO}{\mathcal{W}}  
\newcommand{\SB}{\mathcal{F}}  

\newcommand{\diffOp}{\mathcal{N}}

\newcommand{\dd}{{\bm 1 + \bm e}}

\newcommand{\parSpace}{\Theta}
\newcommand{\parVar}{\bm{\theta}}
\newcommand{\uparVar}{\theta}

\newcommand{\indexset}{\mathcal{I}}
\newcommand{\indexsetsum}{\indexset_{\text{sum}}(w)}
\newcommand{\indexsetsumw}[1]{\indexset_{\text{sum}}(#1)}
\renewcommand{\deg}{d}

\newcommand{\solSpace}{V}
\newcommand{\sol}{u_{\parVar}}
\newcommand{\asol}{\hat u_{\parVar}}
\newcommand{\qoi}{Q}
\newcommand{\aqoi}{\widetilde{\qoi}}
\newcommand{\qoisol}{q}

\newcommand{\sg}{SG\xspace}
\newcommand{\gels}{GELS\xspace}

\newcommand{\poly}{\mathscr{P}}
\newcommand{\polyop}[1]{\mathscr{P}_{\vphantom{\U}\smash{#1}}}
\newcommand{\interpPoints}{\mathcal{T}}
\newcommand{\interpPoint}{\tau}
\newcommand{\lagrange}{\ell}

\DeclareMathOperator{\spn}{span}
\DeclareMathOperator*{\argmin}{argmin}
\newcommand{\norm}[1]{\Vert #1\Vert_\infty}